
\input amstex
\documentstyle{amsppt}

\loadeurm
\loadeusm
\loadbold

\magnification=\magstep1
\hsize=6.5truein
\vsize=9.1truein

\baselineskip=14pt

\def \? {{\bf (???)}}
\def \ideal {\trianglelefteq}
\def \coideal {\dot{\trianglelefteq}}
\def \coleq {\dot{\leq}}
\def \loongrightarrow {\relbar\joinrel\relbar\joinrel\rightarrow}
\def \llongrightarrow
{\relbar\joinrel\relbar\joinrel\relbar\joinrel\rightarrow}
\def \loongtwoheadrightarrow
{\relbar\joinrel\relbar\joinrel\twoheadrightarrow}
\def \llongtwoheadrightarrow
{\relbar\joinrel\relbar\joinrel\relbar\joinrel\twoheadrightarrow}

\def \gerg {\frak g}
\def \gerh {\frak h}
\def \gerk {\frak k}
\def \gerl {\frak l}
\def \gert {\frak t}

\def \gerI {{\frak I}}
\def \gerC {{\frak C}}
\def \calI {{\Cal I}}
\def \calC {{\Cal C}}
\def \gersl {\frak{sl}}
\def \gersln {\frak{sl}_{\,n}}
\def \gergln {\frak{gl}_{\,n}}
\def \gerso {\frak{so}}
\def \h {\hbar}
\def \HA {{\Cal{HA}}}
\def \Hpicc {{\scriptscriptstyle \Cal{H}}}
\def \Ppicc {{\scriptscriptstyle \Cal{P}}}
\def \U {\Bbb U}
\def \F {\Bbb F}
\def \N {\Bbb N}

\def \C {\Bbb C}
\def \kh {\Bbbk [[\h]]}
\def \id {\text{\sl id}}
\def \Ker {\text{\it Ker}}

\def \otimeshat {\,\widehat{\otimes}\,}
\def \otimestilde {\,\widetilde{\otimes}\,}

\def \QUEA {\Cal{Q\hskip1ptUE\hskip-1ptA}}
\def \QFSHA {\Cal{Q\hskip1ptFSHA}}

\def \uhg {U_\h(\gerg)}
\def \fhg {F_\h[[G]]}

\document


\topmatter

\title
  A quantum duality principle for  \\
  coisotropic subgroups and Poisson quotients
\endtitle

\author
       Nicola Ciccoli${}^\dagger$,  \ Fabio Gavarini$^\ddagger$
\endauthor

\leftheadtext{ Nicola Ciccoli, \ \  Fabio Gavarini }
\rightheadtext{ Quantum duality for coisotropic
                subgroups and Poisson quotients }

\affil
  $ {}^\dagger $Universit\`a di Perugia ---
Dipartimento di Matematica  \\
  Via Vanvitelli 1, I-06123 Perugia --- ITALY  \\
                             \\
  \hbox{ $ {}^\ddagger $Universit\`a di Roma ``Tor Vergata'' ---
Dipartimento di Matematica }  \\
  Via della Ricerca Scientifica 1, I-00133 Roma --- ITALY  \\
\endaffil

\address\hskip-\parindent
  Nicola Ciccoli  \newline
  \indent   Universit\`a degli Studi di Perugia   ---   Dipartimento
di Matematica  \newline
  \indent   Via Vanvitelli 1{\,}, I-06123 Perugia --- ITALY{\,},
 \;
   e-mail: \ ciccoli\@{}dipmat.unipg.it  \newline
  {}  \newline
  Fabio Gavarini  \newline
  \indent   Universit\`a degli Studi di Roma ``Tor Vergata''
--- Dipartimento di Matematica  \newline
  \indent   Via della Ricerca Scientifica 1, I-00133 Roma
\!--- ITALY,
 \;
   e-mail: gavarini\@{}mat.uniroma2.it  \newline
\endaddress

\abstract
  We develop a quantum duality principle for coisotropic subgroups of
a (formal) Poisson group and its dual: namely, starting from a quantum
coisotropic subgroup (for a quantization of a given Poisson group)
we provide functorial recipes to produce quantizations of the dual
coisotropic subgroup (in the dual formal Poisson group).  By the natural
link between subgroups and homogeneous spaces, we argue a quantum
duality principle for Poisson homogeneous spaces which are Poisson
quotients, i.e.~have at least one zero-dimensional symplectic leaf.
As an application, we provide an explicit quantization of the
homogeneous  $ {{SL}_n}^{\!\!*} $--space  of Stokes matrices, with
the Poisson structure given by Dubrovin and Ugaglia.
\endabstract

\endtopmatter

%
%
\footnote""{Keywords: \ {\sl Quantum Groups, Poisson Homogeneous
Spaces, Coisotropic Subgroups}.}

\footnote""{ 2000 {\it Mathematics Subject Classification:} \
Primary 17B37, 20G42, 58B32; Secondary 81R50. }

\vskip17pt

\centerline {\bf Introduction }

\vskip11pt

   The natural semiclassical counterpart of the study of quantum
groups is the theory of Poisson groups: indeed, Drinfeld himself
introduced Poisson groups as the semiclassical limits of quantum
groups.  Therefore, it should be no surprise to anyone, anymore,
that the geometry of quantum groups gain in clarity and comprehension
when its connection with Poisson geometry is more transparent.  The
same can be observed when referring to homogeneous spaces.
                                             \par
   In fact, in the study of Poisson homogeneous spaces, a special
r{\^o}le is played by  {\sl Poisson quotients}.  These are those
Poisson homogeneous spaces whose symplectic foliation has at least
one zero-dimensional leaf, so they can be thought of
as pointed Poisson homogeneous spaces, just like Poisson groups
themselves are pointed by the identity element.  When looking at
quantizations of a Poisson homogeneous space, one finds that the
existence is guaranteed only if the space is a quotient (cf.~[EK2]).
Thus the notion of Poisson quotient shows up naturally also from
the point of view of quantization (see [Ci]).
                                    \par
   Poisson quotients are a natural subclass of Poisson homogeneous
$ G $--spaces  ($ G $  a Poisson group),  best adapted to the usual
relation between homogeneous  $ G $--spaces  and subgroups of  $ G \, $:
\, they correspond to  {\sl coisotropic\/}  subgroups.  The quantization
process for a Poisson  $ G $--quotient  then corresponds to a like
procedure for the attached coisotropic subgroup of  $ G $.  Also,
when following an infinitesimal approach one deals with Lie
subalgebras of the Lie algebra  $ \gerg $  of  $ G \, $,  \,
and the coisotropy condition has its natural counterpart in
this Lie algebra setting; the quantization process then is
to be carried on for the Lie subalgebra corresponding to
the initial homogeneous  $ G $--space.
                                    \par
   When quantizing Poisson groups (or Lie bialgebras), a precious tool is
the quantum duality principle (QDP).  Loosely speaking this guarantees
that any quantized enveloping algebra can be turned (roughly speaking)
into a quantum function algebra for the dual Poisson group; viceversa
any quantum function algebra can be turned into a quantization of the
enveloping algebra of the dual Lie bialgebra.  More precisely, let
$ \QUEA $  and  $ \QFSHA $ respectively be the category of all
quantized universal enveloping algebras (QUEA) and the category
of all quantized formal series Hopf algebras (QFSHA), in Drinfeld's
sense.  After its formulation by Drinfeld (see [Dr1], \S 7) the QDP
establishes a category equivalence between  $ \QUEA $  and  $ \QFSHA $
via two functors,  $ \, (\ )' \colon \QUEA \longrightarrow \QFSHA \, $
and  $ \, (\ )^\vee \colon \QFSHA \longrightarrow \QUEA \, $,  \, such
that, starting from a QUEA over a Lie bialgebra (resp.~from a QFSHA over
a Poisson group) the functor  $ (\ )' $  (resp.~$ (\ )^\vee \, $)  gives
a QFSHA (resp.~a QUEA) over the dual Poisson group (resp.~the dual Lie
bialgebra).  In a nutshell,  $ \, {U_\h(\gerg)}' = F_\h[[G^*]] \, $
and  $ \, {F_\h[[G]]}^\vee = U_\h(\gerg^*) \, $  for any Lie bialgebra
$ \gerg \, $.  So from a quantization of any Poisson group this
principle gets out a quantization of the dual Poisson group too.
                                                 \par
   In this paper we establish a similar quantum duality principle for
(closed) coisotropic subgroups of a Poisson group  $ G $,  or equivalently
for Poisson  $ G $--quotients,  sticking to the formal approach which is
best suited for dealing with quantum groups \`a la Drinfeld.  Namely,
given a Poisson group  $ G $  assume quantizations  $ U_\h(\gerg) $  and
$ \fhg $  of it are given; then any formal coisotropic subgroup  $ K $
of  $ G $  has two possible algebraic descriptions via objects related
to  $ U(\gerg) $  or  $ F[[G]] $,  and similarly for the formal
Poisson quotient  $ G\big/K \, $.  Thus the datum of  $ K $  or
equivalently of  $ G\big/K $  is described algebraically in four
possible ways: by quantization of such a datum we mean a quantization
of any one of these four objects.  Our ``QDP'' now is a series of
functorial recipes to produce, out of a quantization of  $ K $  or 
$ G\big/K $  as before, a similar quantization of the so-called 
{\sl complementary dual\/}  of  $ K \, $,  \, i.e.~the coisotropic
subgroup  $ K^\perp $  of  $ G^* $  whose tangent Lie bialgebra
is just  $ \gerk^\perp $  inside  $ \gerg^* \, $,  \, or of the
associated Poisson  $ G^* $--quotient,  namely  $ G^* \big/
K^\perp \, $.
                                                 \par
   We would better stress that, just like the QDP for quantum groups,
ours is by no means an existence result: instead, it can be thought of
as a  {\sl duplication result},  in that it yields a new quantization
(for a complementary dual object) out of one given from scratch.
                                                 \par
   As an aside remark, let us comment on the fact that the more general
problem of quantizing coisotropic manifolds of a given Poisson manifold,
in the context of deformation quantization, has recently raised quite
some interest (see [BGHHW,CF]).
                                                 \par
   As an example, in the last section we show how we can use this quantum
duality principle to derive new quantizations from known ones.  The example
is given by the Poisson structure introduced on the space of Stokes
matrices by Dubrovin (see [Du]) and Ugaglia (see [Ug]) in the framework
of moduli spaces of semisimple Frobenius manifolds.  It was Boalch (cf. [Bo]) that first gave an interpretation of Dubrovin--Ugaglia brackets in terms of Poisson--Lie groups.  We will rather follow later work by Xu (see [Xu])
where it was shown how Boalch construction may be equivalently interpreted
as quotient Poisson structure of the dual Poisson-Lie group  $ G^* $  of the standard  $ SL_n(\Bbbk) $. In more detail the Poisson space of Stokes matrices 
$ \, G^* \big/ H^\perp \, $  is the dual Poisson space to the Poisson space 
$ \, SL_n(\Bbbk) \big/ SO_n(\Bbbk) \, $.
%
%
It has to be noted that the embedding of  $ SO_n(\Bbbk) $  in  $ SL_n(\Bbbk)$  is known to be coisotropic but not Poisson.  Starting, then, from results
obtained by Noumi in [No] related to a quantum version of the embedding
$ \, SO_n(\Bbbk) \lhook\joinrel\longrightarrow SL_n(\Bbbk) \, $  we are
able to interpret them as an explicit quantization of the Dubrovin-Ugaglia
structure.  We provide explicit computations for the case $ \, n = 3
\, $,  \, and draw a sketch with the main guidelines for the general
case.
                                                 \par
   Finally, another, stronger formulation of our QDP for subgroups and
homogeneous spaces can be given in terms of quantum groups of  {\sl
global\/}  type, see [CG].

\vskip1,3truecm

\centerline {\bf \S\; 1 \ The classical setting }

\vskip11pt

   In this section we introduce the notions of Poisson geometry we shall
need in the following: coisotropic subgroups and Poisson quotients, also
called Poisson homogeneous spaces of group type.  Our aim is to stress
their algebraic characterization.

\vskip7pt

  {\bf 1.1 Formal Poisson groups.} \, As already explained, the setup of
the paper is formal geometry.  Recall that a formal variety is uniquely
characterized by a tangent or a cotangent space (at its unique point),
and is described by its ``algebra of regular functions''   --- such
as  $ F[[G]] $  below ---   which is a complete, topological local
ring which can be realized as a  $ \Bbbk $--algebra  of formal power
series.  Hereafter  $ \Bbbk $  is a field of zero characteristic.
                                            \par
     Let  $ \gerg $  be a finite dimensional Lie algebra over  $ \Bbbk
\, $,  and let  $ U(\gerg) $  be its universal enveloping algebra (with
the natural Hopf algebra structure).  We denote by  $ F[[G]] $  the
algebra of functions on the formal algebraic group  $ G $  associated
to  $ \gerg $  (which depends only on  $ \gerg $  itself); this is
a complete, topological Hopf algebra.  One has  $ \, F[[G]] \cong
{U(\gerg)}^* \, $  so that there is a natural pairing of (topological)
Hopf algebras   --- see below ---   between  $ U(\gerg ) $  and
$ F[[G]] \, $.
                                            \par
   In general, if  $ H $,  $ K $  are Hopf algebras (even topological)
over a ring  $ R \, $,  a pairing  $ \, \langle \,\ , \,\ \rangle \,
\colon \, H \times K \longrightarrow R \, $   is called a  {\sl Hopf
pairing\/}  if
  $ \; \big\langle x, y_1 \cdot y_2 \big\rangle = \big\langle \Delta(x),
y_1 \otimes y_2 \big\rangle \, ,  \;\;  \big\langle x_1 \cdot x_2,
y \big\rangle = \big\langle x_1 \otimes x_2, \Delta(y) \big\rangle
\, $,  $ \; \langle x, 1 \rangle = \epsilon(x) \, $,  $ \; \langle 1,
y \rangle = \epsilon(y) \, $, $ \; \big\langle S(x), y \big\rangle =
\big\langle x, S(y) \big\rangle \; $  for all  $ \, x , x_1 , x_2 \in
H \, $,  $ \, y, y_1, y_2 \in K \, $.  Moreover, a pairing is called
{\sl perfect\/}  if it is non-degenerate.  
                                            \par
     Now assume  $ G $  is a formal  {\sl Poisson\/}  (algebraic) group.
Then  $ \gerg $  is a Lie bialgebra,  $ U(\gerg) $  is a co-Poisson Hopf
algebra,  $ F[[G]] $  is a topological Poisson Hopf algebra, and
the Hopf pairing above respects these additional co-Poisson and Poisson
structures.  Furthermore, the linear dual  $ \gerg^* $  of  $ \gerg $
is a Lie bialgebra as well, so a dual formal Poisson group  $ G^* $
exists.

\vskip2pt

   $ \underline{\hbox{\sl Notation}} $:  hereafter, the symbol  $ \,
\coideal \, $  stands for ``coideal'',  $ \, \leq^1 \, $  for ``unital
subalgebra'',  $ \, \coleq \, $  for ``subcoalgebra'',  $ \, \leq_\Ppicc
\, $  for ``Poisson subalgebra'',  $ \, \coideal_\Ppicc \, $  for
``Poisson coideal'',  $ \, \leq_\Hpicc \, $  for ``Hopf subalgebra'',
$ \, \ideal_\Hpicc \, $  for ``Hopf ideal'', and the subscript  $ \ell $
stands for ``left''.  Everything has to be meant in topological sense
if necessary.

\vskip7pt

  {\bf 1.2 Subgroups and homogeneous  $ G $--spaces.}  \, A homogeneous
left  $ G $--space  $ M $  corresponds to a closed subgroup  $ \, K =
K_M \, $,  \, which we assume to be connected, of  $ G $  such that
$ \, M \cong G \big/ K \, $.  Actually, in formal geometry  $ K $  may
be replaced by  $ \, \gerk := \text{\it Lie}(K) \, $  as well.  Then the
whole geometrical setting established by the pair  $ \big( K, G/K \big) $
is algebraically encoded by any one of the following data:

\vskip3pt

   {\it (a)} \  the set  $ \, \calI = \calI(K) \equiv \calI(\gerk) \, $
of all (formal) functions vanishing on  $ K \, $,  that is to say  $ \,
\calI = \big\{ \varphi \!\in\! F[[G]] \,\big|\, \varphi(K) \! = \! 0
\big\} \, $:  \, this is a Hopf ideal of  $ F[[G]] \, $,  in short
$ \; \calI \ideal_\Hpicc \! F[[G]] \; $;
                                                 \par
   {\it (b)} \  the set of all left  $ \gerk $--invariant  functions,
namely  $ \, \calC = \calC(K) \equiv \calC(\gerk) = {F[[G]]}^K \; $:
\, this is a unital subalgebra and left coideal of  $ F[[G]] \, $,
in short  $ \; \calC \leq^1 \! \coideal_\ell \, F[[G]] \; $;
                                                 \par
   {\it (c)} \  the set  $ \, \gerI = \gerI(K) \equiv \gerI(\gerk) \, $
of all left-invariant differential operators on  $ F[[G]] $   which
vanish on  $ \, {F[[G]]}^K \, $,  \, that is  $ \, \gerI = U(\gerg)
\cdot \gerk \, $  (via standard identifications of the set of
left-invariant differential operators with  $ U(\gerg) \, $):
\, this is a left ideal and (two-sided) coideal of  $ U(\gerg)
\, $,  in short  $ \; \gerI(\gerk) = \gerI \ideal_\ell \!
\coideal \, U(\gerg) \, $;
                                                 \par
   {\it (d)} \  the universal enveloping algebra of  $ \, \gerk \, $,
\, denoted  $ \gerC = \gerC(K) \equiv \gerC(\gerk) :=  U(\gerk) \; $:
\, this is a Hopf subalgebra of  $ U(\gerg) \, $,  \, i.e.~$ \; \gerC
\leq_\Hpicc \! U(\gerg) \; $.

\vskip4pt

   In this way any formal subgroup  $ K $  of  $ G \, $,  or the
associated homogeneous  $ G $--space  $ \, G \big/ K \, $,  \, is
characterized   --- via  $ \gerk $  and  $ \gerg $  ---   by any
one of the following algebraic objects:
  $$  \text{\it (a)} \hskip6pt \calI \ideal_\Hpicc \! F[[G]]
\hskip17pt  \text{\it (b)} \hskip6pt \calC \leq^1 \! \coideal_\ell
\, F[[G]]  \hskip17pt  \text{\it (c)} \hskip6pt \gerI \ideal_\ell
\coideal \, U(\gerg)  \hskip17pt  \text{\it (d)} \hskip7pt
\gerC \leq_\Hpicc \! U(\gerg)   \hskip14pt   \eqno (1.1)  $$
   \indent   Clearly  {\it (a)\/}  and  {\it (d)\/}  in (1.1) ideally
focus on the subgroup  $ K \, $,  whereas  {\it (b)\/}  and  {\it (c)\/}
focus more on the formal homogeneous  $ G $--space  $ G \big/ K \, $.
Nevertheless, these four algebraic data are all equivalent to each
other.  To express this algebraically, we need some more notation.
                                                 \par
   For any Hopf algebra  $ H \, $,  with counit  $ \epsilon \, $,  and
every submodule  $ \, M \subseteq H \, $,  \, we set:  $ \, M^+ := M
\cap \text{\it Ker}\,(\epsilon) \, $  and  $ \, H^{\text{\it co}M} :=
\big\{\, y \in H \,\big|\, \big( \Delta(y) - y \otimes 1 \big) \in H
\otimes M \,\big\} \, $  (the set of  $ M $--{\sl coinvariants\/}  of
$ H \, $).  Letting  $ \Bbb{A} $  be the set of all subalgebras left
coideals of  $ H $  and  $ \Bbb{K} $  be the set of all coideals
left ideals of  $ H \, $,  we have well-defined maps  $ \, \Bbb{A}
\longrightarrow \Bbb{K} \, $,  $ \, A \mapsto H \cdot A^+ \, $,
and  $ \, \Bbb{K} \longrightarrow \Bbb{A} \, $,  $ \, K \mapsto
H^{\text{\it co}K} \, $  (cf.~[Ma], and references therein).
                                                 \par
   Then the above mentioned equivalence stems from the following
relations, which starting from any one of the four items in (1.1)
allow one to reconstruct the remaining ones:

 \vskip2pt

   {\it --- (1)}  \,  {\sl orthogonality relations\/}   --- w.r.t.~the
natural pairing between  $ F[[G]] $  and  $ U(\gerg) $  ---   namely
$ \, \calI = \gerC^\perp $,  $ \, \gerC = \calI^\perp $,  \, linking
{\it (a)\/}  and  {\it (d)},  and  $ \, \calC = \gerI^\perp $,  $ \,
\gerI = \calC^\perp $,  linking  {\it (b)\/}  and  {\it (c)\/};

 \vskip2pt

   {\it --- (2)}  \,  {\sl subgroup-space correspondence},
namely  $ \, \calI = F[[G]] \cdot \calC^+ \, $,  $ \, \calC
= {F[[G]]}^{\text{\it co}\calI} \, $,  \, linking  {\it (a)\/}  and
{\it (b)},  and  $ \; \gerI = U(\gerg) \cdot \gerC^+ \, $,  $ \, \gerC
= {U(\gerg)}^{\text{\it co} \gerI} \, $,  \, linking  {\it (c)\/}  and
{\it (d)}.  Moreover, the maps  $ \, \Bbb{A} \longrightarrow \Bbb{K} \, $
and  $ \, \Bbb{K} \longrightarrow \Bbb{A} \, $  considered above are
inverse to each other in the formal setting.

 \vskip7pt

  {\bf 1.3 Coisotropic subgroups and Poisson quotients.}  \, When
$ G $  is a Poisson group, a distinguished class of subgroups   ---
the  {\sl coisotropic\/}  ones  ---   is of special interest.
                                                 \par
   A closed formal subgroup  $ K $  of  $ G $  with Lie algebra  $ \gerk $
is called  {\sl coisotropic}  if its defining ideal  $ {\Cal I}(\gerk) $
is a (topological) Poisson subalgebra of  $ F[[G]] \, $.  The following
are equivalent:

\vskip3pt

   {\it (C-i)}  \,  $ K $  is a coisotropic formal subgroup of  $ G \, $;
                                                 \par
   {\it (C-ii)}  \,  $ \delta(\gerk) \subseteq \gerk \wedge \gerg \, $,
that is  $ \gerk $  is a Lie coideal of  $ \gerg \, $;
                                                 \par
   {\it (C-iii)}  \,  $ \gerk^\perp \, $  is a Lie subalgebra
of  $ \gerg^* \, $

\vskip3pt

\noindent   (see [Lu]).  Clearly  {\it (C-ii)\/}  and  {\it (C-iii)\/}
characterize coisotropic subgroups in algebraic terms.
                                                 \par
   As for homogeneous spaces, recall that a formal Poisson manifold  $ \,
(M, \omega_M) \, $  is a  {\sl Poisson homogeneous  $ G $--space\/}  if
there is a smooth homogeneous action  $ \, \phi \colon \, G \times M
\rightarrow M \, $  which is a Poisson map with respect to the product
Poisson structure.
                                                 \par
   In addition,  $ \, (M, \omega_M) \, $  is said to be  {\sl of group
type\/}  (after Drinfeld [Dr2]), or simply a  {\sl Poisson quotient},
if there exists a coisotropic closed Lie subgroup  $ K_M $  of  $ \, G $
such that  $ \, G \big/ K_M \simeq M \, $  and the natural projection
$ \, \pi \colon \, G \longrightarrow G \big/ K_M \simeq M \, $  is a
Poisson map.
                                                 \par
   The following is a characterization of Poisson quotients (cf.~[Za]):

\vskip3pt

   {\it (PQ-i)}  \,\; there exists  $ \, x_0 \in M \, $  such that its
stabilizer  $ \, G_{x_0} \, $  is coisotropic in  $ \, G \, $;
                                                 \par
   {\it (PQ-ii)}  \,\, there exists  $ \, x_0 \in M \, $  such that
$ \, \phi_{x_0} \, \colon \, G \longrightarrow M \, $,  $ \, g \mapsto
\phi(g,x_0) \, $,  \; is a Poisson map, that is  $ M $  is a Poisson
quotient;
                                                 \par
   {\it (PQ-iii)}  \, there exists  $ \, x_0 \in M \, $  such that
$ \, \omega_M(x_0) = 0 \, $.

\vskip4pt

  {\it  $ \underline{\hbox{\it Remark}} $:} \,  in Poisson geometry,
the usual relationship between closed subgroups of  $ G $  and
$ G $--homogeneous spaces does not hold anymore.  In fact, in the
{\sl same\/}  conjugacy class one can have Poisson subgroups,
coisotropic subgroups  {\sl and\/}  non-coisotropic subgroups.
We saw above that Poisson quotients correspond to Poisson
homogeneous spaces in which at least one of the stabilizers is
coisotropic; many such examples can be found, for instance, in [LW].
On the other hand many interesting Poisson homogeneous spaces are not
of group type, as it is the case for covariant (in particular invariant)
symplectic structures.   $ \quad \diamondsuit $

\vskip7pt

\proclaim{Definition 1.4}
                                          \hfill\break
   \indent   (a) \, If  $ K $  is a formal coisotropic subgroup of
$ \, G $,  we call  {\sl complementary dual}  of  $ \, K $  the formal
subgroup  $ K^\perp $  of  $ \, G^* $  whose tangent Lie algebra is
$ \, \gerk^\perp $  (with  $ \, G^* $  as in \S 1.1).
                                          \hfill\break
   \indent   (b) \, If  $ \, M \cong G \big/ K_M \, $  is a formal
Poisson  $ G $--quotient,  with  $ K_M $  coisotropic, we call  {\sl
complementary dual}  of  $ M $  the formal Poisson  $ G^* $--quotient
$ \, M^\perp := G^* \big/ K_{\!M}^{\;\perp} \, $.
\endproclaim

\vskip5pt

  {\bf 1.5 Remarks:} \, {\it (a)} \, The fact to be highlighted in
the above definition is that a subset  $ \gerk $  of  $ \gerg $  is
a Lie coideal if and only if  $ \gerk^\perp $  is a Lie subalgebra
of  $ \gerg^* \, $.  This is why we have dual Poisson quotients.
Even more, by  {\it (C-i,ii,iii)\/}  in \S 1.3, the complementary
dual subgroup to a coisotropic subgroup is  {\sl coisotropic\/}  too,
and taking twice the complementary dual gives back the initial subgroup.
Similarly, the Poisson homogeneous space which is complementary dual
to a Poisson homogeneous space of group type is in turn  {\sl of group
type\/}  as well, and taking twice the complementary dual gives back
the initial manifold.  So Definition 1.4 makes sense, and the notion
of complementary duality is self-dual, in both cases.
                                            \par
   {\it (c)} \, The notion of Poisson homogeneous  $ G $--spaces of group
type was first introduced by Drinfeld in [Dr2]: here the relation between
such  $ G $--spaces  and Lagrangian subalgebras of Drinfeld's double
$ \, D(\gerg) = \gerg \oplus \gerg^* \, $  is also explained.  This
is further developed in [EL].
                                            \par
   {\it (d)} \, We denote by  $ \, \hbox{\it
co}\hskip1pt\Cal{S}(G) \, $  the set of all formal coisotropic subgroups
of  $ G \, $,  \, which is as well described by the set of all Lie
subalgebras, Lie coideals of  $ \gerg \, $.  This is a lattice
{w.~r.~t.}  set-theoretical inclusion, hence it can (and will)
also be thought of as a category.   $ \quad \diamondsuit $

\vskip6pt

  {\bf 1.6 Algebraic characterization of coisotropic subgroups.} \,
Let  $ K $  be a formal co\-isotropic subgroup of  $ G $.  Taking
$ \gerI $,  $ \calC $,  $ \calI $  and  $ \gerC $  as in \S 1.2,
coisotropy corresponds to
 \vskip-16pt
  $$  \text{\it (a)} \;\, \calI \leq_\Ppicc F[[G]] \, ,  \hskip21pt
\text{\it (b)} \;\, \calC \leq_\Ppicc F[[G]] \, ,  \hskip21pt
\text{\it (c)} \;\, \gerI \, \coideal_\Ppicc \, U(\gerg) \, ,
\hskip21pt  \text{\it (d)} \;\, \gerC \, \coideal_\Ppicc \, U(\gerg)  $$
 \vskip-4pt
Thus a formal coisotro\-pic subgroup of  $ G $  is  identified by any
one of the algebraic objects
 \vskip-13pt
  $$  \text{\it (a)} \; \calI \ideal_\Hpicc \leq_\Ppicc \! F[[G]] \, ,
\hskip3pt  \text{\it (b)} \; \calC \leq^1 \! \coideal_\ell \leq_\Ppicc
\! F[[G]] \, ,  \hskip3pt  \text{\it (c)} \; \gerI \ideal_\ell \!
\coideal \, \coideal_\Ppicc \, U(\gerg) \, ,  \hskip3pt  \text{\it (d)}
\; \gerC \leq_\Hpicc \! \coideal_\Ppicc \, U(\gerg) \, .
\hskip5pt (1.2)  $$
 \vskip-1pt
   \indent   Note also that  $ K $  being coisotropic reflects the fact
that the distinguished point  $ eK $ (where  $ \, e \in G \, $  is the
identity element) in the formal Poisson  $ G $--space  $ G\big/K $  is
a zero-dimensional leaf.  Then the algebra of regular functions on
$ G\big/K \, $,  \, already realized as  $ {F[[G]]}^K $,  will be also
denoted by  $ F \big[\big[ G\big/K \big]\big] \, $.  Moreover, we can
always choose a system of parameters for  $ G \, $,  say  $ \, \big\{
j_1, \dots, j_k, j_{k+1}, \dots, j_n \big\} \, $  such that  $ \, k =
\dim(K) \, $,  $ \, n = \dim(G) \, $,  $ \, {F[[G]]}^K \! = \Bbbk[[
j_{k+1},\dots,j_n]] $  (the topological  {\sl subalgebra\/}  of
$ F[[G]] $  generated by  $ \big\{ j_{k+1}, \dots, j_n \big\} \! $)
and  $ \calI(K) = \! \big( j_{k+1}, \dots, j_n \big) $  (the
{\sl ideal\/}  of  $ F[[G]] $  generated by  $ \big\{ j_{k+1},
\dots, j_n \big\} $).

\vskip1,1truecm

\centerline {\bf \S\; 2 \ The quantum setting }

\vskip10pt

   This section is devoted to recall quantum groups and Drinfeld's
QDP for quantum groups, to introduce our concept of quantization for
coisotropic subgroups and Poisson quotients, and to explain the basic
idea of our QDP for the latters.

\vskip6pt

  {\bf 2.1 Topological  $ \Bbbk[[\h]] $--modules  and tensor structures.}
\, Let  $ \kh $  be the topological ring of formal power series in the
indeterminate  $ \h \, $.  If  $ X $  is any  $ \kh $--module, we
set  $ \, X_0 := X \big/ \h X = \Bbbk \otimes_{\kh} X \, $, \, the
{\sl specialization\/}  of  $ X $  at  $ \, \h = 0 \, $,  or
{\sl semiclassical limit\/}  of  $ X \, $.
                                            \par
   Let  $ {\Cal T}_{\otimeshat} $  be the category whose objects are all
topological  $ \kh $--modules  which are topologically free and whose
morphisms are the  $ \kh $--linear  maps (which are automatically
continuous).  It is a tensor category for the
tensor product  $ \, T_1 \otimeshat T_2 \, $  defined as the separated
$ \h $--adic  completion of the algebraic tensor product  $ \, T_1
\otimes_{\kh} T_2 \, $  (for all  $ T_1 $,  $ T_2 \in
{\Cal T}_{\otimeshat} $).  We denote by  $ \HA_{\otimeshat} $
the subcategory of  $ {\Cal T}_{\otimeshat} $  whose objects are all
the Hopf algebras in  $ {\Cal T}_{\otimeshat} $  and whose morphisms
are all the Hopf algebra morphisms in $ {\Cal T}_{\otimeshat} $.
                                            \par
   Let  $ {\Cal P}_{\otimestilde} $  be the category whose objects are
all topological  $ \kh $--modules  isomorphic to modules of the type
$ {\kh}^E $  (with the Tikhonov product topology) for some set  $ E \, $,
\, and whose morphisms are the  $ \kh $--linear  continuous maps.  It
is a tensor category w.r.t.~the tensor product  $ \, P_1 \otimestilde
P_2 \, $  defined as the completion of the algebraic tensor product
$ \, P_1 \otimes_{\kh} P_2 \, $  w.r.t.~the weak topology:  thus
$ \, P_i \cong {\kh}^{E_i} $  ($ i = 1 $,  $ 2 $)  yields  $ \, P_1
\otimestilde P_2 \cong {\kh}^{E_1 \times E_2} \, $  (for all  $ P_1 $,
$ P_2 \in {\Cal P}_{\otimestilde} $).  We call  $ \HA_{\otimestilde} $
the subcategory of  $ {\Cal P}_{\otimestilde} $  whose objects are all
the Hopf algebras in  $ {\Cal P}_{\otimestilde} $  and whose morphisms
are all the Hopf algebra morphisms in  $ {\Cal P}_{\otimestilde} $.

\vskip6pt

\proclaim{Definition 2.2}  (cf.~[Dr1, \S~7])
                                         \hfill\break
  \indent  (a) \, We call QUEA any  $ \, H \in \HA_{\otimeshat} $
such that  $ \, H_0 := H \big/ \h H \, $  is a co-Poisson Hopf algebra
isomorphic to  $ U(\gerg) $  for some finite dimensional Lie bialgebra
$ \gerg $  (over  $ \Bbbk $);  in this case we write  $ \, H =
U_\h(\gerg) \, $,  \, and say  $ H $  is a  {\sl quantization\/}
of  $ U(\gerg) $.  We call $ \QUEA $  the full tensor subcategory
of  $ \, \HA_{\otimeshat} $  whose objects are QUEA, relative to
all possible  $ \gerg $  (see also Remark 2.3 below).
                                         \hfill\break
  \indent  (b) \, We call QFSHA any  $ \, K \in \HA_{\otimestilde} $ 
such that  $ \, K_0 := K \big/ \h K \, $  is a topological Poisson
Hopf algebra isomorphic to  $ F[[G]] $  for some finite dimensional
formal Poisson group  $ G $  (over  $ \Bbbk $);  then we write  $ \,
H = F_\h[[G]] \, $,  \, and say  $ K $  is a  {\sl quantization\/} 
of  $ F[[G]] $.  We call  $ \QFSHA $  the full tensor subcategory
of  $ \, \HA_{\otimestilde} $  whose objects are QFSHA, relative
to all possible  $ G $  (see also Remark 2.3 below).
\endproclaim

\vskip1pt

  {\bf Remarks 2.3:}  \, If  $ \, H \in \HA_{\otimeshat} $
is such that  $ \, H_0 := H \big/ \h H \, $  as a Hopf algebra
is isomorphic
 to  $ U(\gerg) $  for some Lie algebra  $ \gerg $,
then  $ \, H_0 = U(\gerg) \, $  is also a  {\sl co-Poisson\/}
Hopf algebra w.r.t.~the
 \eject
\noindent
 Poisson cobracket  $ \delta $  defined as follows: if  $ \,
x \in H_0 \, $  and  $ \, x' \in H \, $  gives  $ \, x = x' + \h \, H
\, $,  \, then  $ \, \delta(x) := \big( \h^{-1} \, \big( \Delta(x') -
\Delta^{\text{op}}(x') \big) \big) + \h \, H \otimeshat H \, $;  \,
then (by [Dr1, \S 3, Theorem 2]) the restriction of  $ \delta $
makes  $ \gerg $  into a Lie bialgebra.  Similarly, if  $ \, K \in
\HA_{\otimestilde} $  is such that  $ \, K_0 := K \big/ \h K \, $ 
is a topological Poisson Hopf algebra isomorphic to  $ F[[G]] $  for
some formal group  $ G $  then  $ \, K_0 = F[[G]] \, $  is also a
topological  {\sl Poisson\/}  Hopf algebra w.r.t.~the Poisson bracket
$ \{\,\ ,\ \} $  defined as follows: if  $ \, x $,  $ y \in K_0 \, $
and  $ \, x' $,  $ y' \in K \, $  give  $ \, x = x' + \h \, K $,  $ \,
y = y' + \h \, K $,  \, then  $ \, \{x,y\} := \big( \h^{-1} (x' \, y' -
y' \, x') \big) + \h \, K \, $;  \, then  $ F[[G]] $  is (the algebra
of regular functions on) a  {\sl Poisson\/}  formal group.  These
natural co-Poisson and Poisson structures are the ones considered
in Definition 2.2 above.

\vskip7pt

  {\bf 2.4 Drinfeld's functors.} \,  Let  $ H $  be a (topological)
Hopf algebra over  $ \kh $.  For each  $ \, n \in \N \, $,  define
$ \; \Delta^n \colon H \longrightarrow H^{\otimes n} \; $  by  $ \,
\Delta^0 := \epsilon \, $,  $ \, \Delta^1 := \id_{\scriptscriptstyle H} $,
\, and  $ \, \Delta^n := \big( \Delta \otimes \id_{\scriptscriptstyle
H}^{\,\otimes (n-2)} \big) \circ \Delta^{n-1} \, $  if  $ \, n \geq 2
\, $.  For any ordered subset  $ \, E = \{i_1, \dots, i_k\} \subseteq
\{1, \dots, n\} \, $  with  $ \, i_1 < \dots < i_k \, $,  \, define the
morphism  $ \; j_{\scriptscriptstyle E} : H^{\otimes k} \longrightarrow
H^{\otimes n} \; $  by  $ \; j_{\scriptscriptstyle E} (a_1 \otimes
\cdots \otimes a_k) := b_1 \otimes \cdots \otimes b_n \; $  with
$ \, b_i := 1 \, $  if  $ \, i \notin \Sigma \, $  and  $ \, b_{i_m}
:= a_m \, $  for  $ \, 1 \leq m \leq k \, $;  then set  $ \; \Delta_E
:= j_{\scriptscriptstyle E} \circ \Delta^k \, $,  $ \; \Delta_\emptyset
:= \Delta^0 \; $  and  $ \; \delta_E := {\textstyle \sum\limits_{E'
\subset E}} {(-1)}^{n- \left| E' \right|} \Delta_{E'} \, $,  $ \;
\delta_\emptyset := \epsilon \; $.  The inverse formula  $ \; \Delta_E
= \sum_{\Psi \subseteq E} \delta_\Psi \, $  holds too.  We shall also
use the notation  $ \, \delta_0 := \delta_\emptyset \, $,  $ \, \delta_n
:= \delta_{\{1, 2, \dots, n\}} \, $.  Then we define
  $$  H' := \big\{\, a \in H \,\big\vert\; \delta_n(a) \in
h^n H^{\otimes n} \; \forall\, n \in \N \,\big\}  \qquad
\big( \subseteq H \, \big ) \; .   $$
   \indent   Note that the useful formula  $ \; \delta_n =
{(\id_{\scriptscriptstyle H} - \epsilon)}^{\otimes n} \circ
\Delta^n \; $  holds, for all  $ \, n \in \N_+ \, $.  Since  $ H $
splits as  $ \, H = \kh \cdot 1_{\scriptscriptstyle H} \oplus
J_{\scriptscriptstyle H} \, $,  and  $ \, (\text{\sl id} -
\epsilon) \, $  projects  $ H $  onto  $ \, J_{\scriptscriptstyle H}
:= \text{\it Ker}\,(\epsilon) \, $,  from  $ \; \delta_n =
{(\id_{\scriptscriptstyle H} - \epsilon)}^{\otimes n} \circ \Delta^n
\; $  we get  $ \; \delta_n(a) = {(\id_{\scriptscriptstyle H} -
\epsilon)}^{\otimes n} \big( \Delta^n(a) \big) \in
{J_{\scriptscriptstyle H}}^{\!\otimes n} \; $
for all  $ \, a \in H \, , \, n \in \N \, $.
                                        \par
   For later use, we recall that ([KT, Lemma 3.2]), if  $ \Phi $  is
any finite subset of  $ \N $  then
  $$  \displaylines{
   \hfill   \delta_\Phi(ab) \; = \; {\textstyle \sum_{\Lambda
\cup Y = \Phi}} \hskip2pt \delta_\Lambda(a) \, \delta_Y(b)
\qquad  \forall\;\; a, b \in H \; ;   \hfill (2.1)  \cr
   \hfill   \delta_\Phi(ab - ba) \; = \;
{\textstyle \sum_{\Sb \Lambda \cup Y = \Phi  \\
\Lambda \cap Y \not= \emptyset  \endSb}} \hskip1pt
\big( \delta_\Lambda(a) \, \delta_Y(b) - \delta_Y(b) \,
\delta_\Lambda(a) \big) \quad  \forall\; a, b \in H \, ,
\; \Phi \not= \emptyset \; .
\hfill (2.2)  \cr }  $$
   \indent   Now let  $ \, I_{\scriptscriptstyle H} := \epsilon^{-1}
\big( \h \, \kh \big) \, $;  \, set  $ \, H^\times :=
\sum\limits_{n \geq 0} \h^{-n} {I_{\scriptscriptstyle H}}^{\!n}
= \sum\limits_{n \geq 0} {\big( \h^{-1} I_{\scriptscriptstyle H}
\big)}^n = \bigcup\limits_{n \geq 0} {\big( \h^{-1}
I_{\scriptscriptstyle H} \big)}^n = \sum_{n \geq 0} \h^{-n}
{J_{\scriptscriptstyle H}}^{\!n} \, $  (inside  $ \, \Bbbk((\h))
\otimes_{\kh} H \, $),  and define
  $$  H^\vee :=  \h\text{--adic completion of the
$ \kh $--module }  H^\times  \; .  $$
   \indent   By means of this constructions, the QDP says that any
QUEA provides also a QFSHA for the dual Poisson group, and any QFSHA
yields also a QUEA for the dual Lie bialgebra:

\vskip7pt

\proclaim {Theorem 2.5} {\sl (``The quantum duality principle'' [=QDP];
cf.~Drinfel'd [Dr1, \S 7]; see also Etingof and Schiffman [ES, \S 10.2],
or Gavarini [Ga1], for a proof)}
\, The assignments  $ \, H \mapsto H^\vee \, $  and
$ \, H \mapsto H' \, $,  respectively, define tensor functors
$ \, \QFSHA \longrightarrow \QUEA \, $  and  $ \, \QUEA
\longrightarrow \QFSHA \, $,  \, which are inverse to each
other.  Indeed,  for all  $ \, U_\h(\gerg) \in \QUEA \, $  and all
$ \, F_\h[[G]] \in \QFSHA \, $  one has
  $$  {U_\h(\gerg)}' \Big/ \h \, {U_\h(\gerg)}' = F[[G^*]] \, ,
\qquad  {F_\h[[G]]}^\vee \Big/ \h \, {F_\h[[G]]}^\vee =
U(\gerg^*)  $$
that is, if  $ \, U_\h(\gerg) $  is a quantization of  $ \, U(\gerg) $
then  $ \, {U_\h(\gerg)}' $  is a quantization of  $ \, F[[G^*]] $,
and if  $ \, F_\h[[G]] $  is a quantization of  $ \, F[[G]] $
then  $ \, {F[[G^*]]}^\vee $  is a quantization of
$ \, U(\gerg^*) \, $.   \qed
\endproclaim

\vskip7pt

   In addition, Drinfeld's functors respect Hopf duality, in the sense
of the following
 \eject
%
%

\proclaim {Proposition 2.6} {\sl (see Gavarini [Ga1, Proposition 2.2])}
\, Let  $ \, U_\h \in \QUEA \, $,  $ \, F_\h \in \QFSHA $  and let
$ \, \pi \, \colon \, U_\h \times F_\h \rightarrow \kh \, $  be a
perfect Hopf pairing whose specialization at  $ \, \h = 0 \, $  is
perfect as well.  Then  $ \pi $  induces   --- by restriction on
l.h.s.~and scalar extension on r.h.s.~---   a perfect Hopf pairing
$ \, {U_\h}' \times {F_\h}^{\!\vee} \!\rightarrow \kh \, $  whose
specialization 
   \hbox{at  $ \h = 0 $  is again perfect too.   $ \! \square $}   
\endproclaim

\vskip6pt

  {\bf 2.7 Quantum subgroups and quantum homogeneous spaces.} \, From
now on, let  $ G $  be a formal Poisson group,  $ \, \gerg := \text{\it
Lie}(G) $  its tangent Lie bialgebra.  We assume a quantization of  $ G $
is given, in the sense that a QFSHA  $ \fhg $  quantizing  $ F[[G]] $
and a QUEA  $ \uhg $  quantizing  $ U(\gerg) $  are given such that, in
addition,  $ \, \fhg \cong {\uhg}^* := \text{\sl Hom}_{\,\kh}\big(\uhg,
\kh\big) \, $  as topological Hopf algebras; the latter requirement is
equivalent to fix a perfect Hopf algebra pairing between  $ \fhg $  and
$ \uhg \, $  whose specialization at  $ \, \h = 0 \, $  be perfect too.
Note that this assumption is not restrictive: by [EK1], a QUEA  $ \uhg $
as required always exists, and then  $ \fhg $  can be simply taken to
be  $ \, \fhg \cong {\uhg}^* \, $,  \, by definition.  Finally, as
a matter of notation we denote by  $ \, \pi_{F_\h} \colon \fhg
\loongtwoheadrightarrow F[[G]] \, $   and  $ \, \pi_{U_\h} \colon
\uhg \loongtwoheadrightarrow U(\gerg) \, $  the specialization maps,
and we set  $ \, F_\h := \fhg \, $,  $ \, U_\h := \uhg \, $.
                                               \par
   Let  $ K $  be a formal subgroup of  $ G \, $,  \, and  $ \, \gerk :=
\text{\it Lie}(K) \, $.  As quantization of  $ K $  and/or of  $ \, G
\big/ K \, $,  \, we mean a quantization of any one of the four algebraic
objects  $ \calI $,  $ \calC $,  $ \gerI $  and  $ \gerC $  associated
to them in \S 1.2, that is either of the following:
 \vskip-11pt
  $$  \hbox{ \hskip-3pt   $ \eqalign{
   (a)  &  \, \text{\ a left ideal, coideal \ }
\calI_\h \ideal_\ell \coideal \; F_\h[[G]]  \, \text{\ such that \ }
\calI_\h \big/ \h \, \calI_\h \, \cong \, \pi_{F_\h}(\calI_\h)
\, = \, \calI  \cr
   (b)  &  \, \text{\ a subalgebra, left coideal \ }
\calC_\h \leq^1 \coideal_\ell \, F_\h[[G]]  \, \text{\ such that \ }
\calC_\h \big/ \h \, \calC_\h \, \cong \, \pi_{F_\h}(\calC_\h)
\, = \, \calC  \cr
   (c)  &  \, \text{\ a left ideal, coideal \ }
\gerI_\h \ideal_\ell \coideal \;\; U_\h(\gerg)  \, \text{\ such that\ } \,
\gerI_\h \big/ \h \, \gerI_\h \, \cong \, \pi_{U_\h}(\gerI_\h)
\, = \, \gerI  \cr
   (d)  &  \, \text{\ a subalgebra, left coideal \ }
\gerC_\h \leq^1 \coideal_\ell \, U_\h(\gerg)  \, \text{\ such that\ } \,
\gerC_\h \big/ \h \, \gerC_\h \, \cong \, \pi_{U_\h}(\gerC_\h)
\, = \, \gerC  \cr } $ }   \hskip0pt (2.3)  $$
 \vskip-3pt
In (2.3) the constraint  $ \; \calI_\h \big/ \h \, \calI_\h \, \cong
\, \pi_{F_\h}(\calI_\h) \, = \, \calI \; $  means the following.  By
construction  $ \; \calI_\h \! \lhook\joinrel\loongrightarrow \! \fhg
\,{\buildrel {\pi_{F_\h}} \over \llongtwoheadrightarrow}\, \fhg \Big/
\h \, \fhg \, \cong \, F[[G]] \, $,  \, and the composed map
   \hbox{$ \, \calI_\h \! \loongrightarrow \! F[[G]] $}
 factors through  $ \, \calI_\h
\big/ \h \, \calI_\h \, $;  \, then we ask that the induced map  $ \,
\calI_\h \big/ \h \, \calI_\h \loongrightarrow F[[G]] \, $  be a
bijection onto  $ \pi_{F_\h}(\calI_\h) \, $,  and that the latter do
coincide with  $ \calI \, $;  \; of course this bijection will also
respects all Hopf operations, because  $ \pi_{F_\h} $  does.  Similarly
for the other conditions.
                                            \par
   The existence of any of such objects is a separate problem, which
we shall not tackle.  However, the four existence problems are in fact
equivalent, in that as one solves any one of them, a solution follows
for the remaining ones.  Indeed, much like in \S 1.2, one has:

 \vskip3pt

   {\it ---  $ \underline{\hbox{{\it (a)}  $ \Longleftrightarrow $
{\it (d)}  \, and \,  {\it (b)}  $ \Longleftrightarrow $  {\it (c)}}} $:}
\, if  $ \calI_\h $  exists as in  {\it (a)},  then  $ \, \gerC_\h :=
{\calI_\h}^{\!\perp} \, $  enjoys the properties in  {\it (d)\/};  \,
conversely, if  $ \gerC_\h $  exists as in  {\it (d)},  then  $ \,
\calI_\h := {\gerC_\h}^{\!\perp} \, $  enjoys the properties in
{\it (a)\/}  (hereafter orthogonality is meant w.r.t.~the fixed
Hopf pairing between  $ \fhg $  and  $ U_\h(\gerg) \, $).  The
equivalence \,  {\it (b)}  $ \Longleftrightarrow $  {\it (c)}
\; follows from a like orthogonality argument.

 \vskip1pt

   {\it ---  $ \underline{\hbox{{\it (a)}  $ \Longleftrightarrow $
{\it (b)}  \, and \,  {\it (c)}  $ \Longleftrightarrow $  {\it (d)}}} $:}
\, if  $ \calI_\h $  exists as in  {\it (a)},  then  $ \, \calC_\h :=
\calI_\h^{\,\text{\it co}\calI_\h} \, $  is an object like in  {\it
(b)\/};  \, on the other hand, if  $ \calC_\h $  as in  {\it (b)\/}
is given, then  $ \, \calI_\h := \fhg \cdot \calC_\h^{\,+} \, $
enjoys all properties in  {\it (a)\/}  (notation of \S 1.2).  The
equivalence \,  {\it (c)}  $ \! \Longleftrightarrow \! $  {\it (d)}
   \hbox{\; stems from a like argument.}

\vskip2pt

   From now on,  {\sl we assume from scratch that quantizations
$ \calI_\h \, $,  $ \calC_\h \, $,  $ \gerI_\h $  and  $ \gerI_\h $
as in (2.3) be given, and that they be linked by the like of relations
{\it (1)--(2)\/}  in \S 1.2, namely}
  $$  \hbox{ $ \matrix
   \hbox{\it (i)}  \;\quad  \calI_\h = {\gerC_\h
\phantom{)}}^{\!\!\!\!\perp} \, ,  \;\quad  \gerC_\h
= {\calI_\h \phantom{)}}^{\!\!\!\perp}  \qquad
  &  \hbox{\it (ii)}  \;\quad  \gerI_\h = {\calC_\h
\phantom{)}}^{\!\!\!\perp} \, ,  \;\quad  \calC_\h =
{\gerI_\h \phantom{)}}^{\!\!\!\!\perp}  \\
   \hbox{\it (iii)}  \quad  \calI_\h = F_\h \cdot \calC_\h^{\,+} \; ,
\;\;  \calC_\h = {F_\h}^{\!\text{\it co}\calI_\h}  \quad
  &  \hbox{\it (iv)}  \quad  \gerI_\h = U_\h \cdot \gerC_\h^{\,+} \; ,
\;\;  \gerC_\h = {U_\h}^{\!\text{\it co}\gerI_\h}  \\
   \endmatrix $ }   \eqno (2.4)  $$
   \indent   In fact, one of the objects is enough to have all the
others, in such a way that the previous assumption holds.  Indeed, if
$ \, \hbox{\it co}\hskip1pt\Cal{S} := \hbox{\it co}\hskip1pt\Cal{S}(G)
\, $  let  $ \, Y_\h \big( \hbox{\it co}\hskip1pt\Cal{S} \big) :=
\big\{ Y_\h(\gerk) \big\}_{\gerk \in \, co\Cal{S}} \, $  for all
$ \, Y \in \big\{ \calI, \calC, \gerI, \gerC \big\} \, $.  The
equivalences  {\it (a)}  $ \Longleftrightarrow $  {\it (d)},
{\it (b)}  $ \Longleftrightarrow $  {\it (c)},  {\it (a)}
$ \Longleftrightarrow $  {\it (b)}  \, and  {\it (c)}
$ \Longleftrightarrow $  {\it (d)}  \, seen above are
given by bijective maps  $ \, \calI_\h \big( \hbox{\it
co}\hskip1pt\Cal{S} \big) \longleftrightarrow \gerC_\h
\big( \hbox{\it co} \hskip1pt\Cal{S} \big) \, $,  $ \,
\calC_\h \big( \hbox{\it co}\hskip1pt\Cal{S} \big)
\longleftrightarrow \gerI_\h \big( \hbox{\it co}
\hskip1pt\Cal{S} \big) \, $,  $ \, \calI_\h \big(
\hbox{\it co}\hskip1pt\Cal{S} \big) \longleftrightarrow
\calC_\h \big( \hbox{\it co}\hskip1pt\Cal{S} \big) \, $
and  $ \, \gerI_\h \big( \hbox{\it co}\hskip1pt\Cal{S}
\big) \longleftrightarrow \gerC_\h \big( \hbox{\it co}
\hskip1pt\Cal{S} \big) \, $  respectively.  Altogether these
maps form a square, which happens to be  {\sl commutative}.
This follows from the fact that each of these maps, or their
inverse, is of type  $ \, X_\h \mapsto X_\h^\perp \, $,  $ \, A_\h
\mapsto H_\h A_\h^+ \, $  or  $ \, K_\h \mapsto H_\h^{\text{\it co}K_\h}
\, $  (see \S 1.2): since the general relations  $ \, X_\h \subseteq
\big(X_\h^\perp\big)^\perp \, $  and  $ \, A_\h \subseteq H_\h^{\text{\it
co}(H_\h A_\h^+)} \, $  hold, and these inclusions turn to identities at
$ \, \h = 0 \, $,  one gets  $ \, X_\h = \big( X_\h^\perp \big)^\perp \, $
and  $ \, A_\h = H_\h^{\text{\it co}(H_\h A_\h^+)} \, $,  which are the
key steps to prove (easily) that the square of maps is commutative, as
claimed.
                                          \par
   Note also that the sets  $ \calI_\h \big( \hbox{\it co}\hskip1pt
\Cal{S} \big) $,  $ \gerC_\h \big( \hbox{\it co}\hskip1pt\Cal{S} \big) $,
$ \calC_\h \big( \hbox{\it co}\hskip1pt\Cal{S} \big) $  and  $ \gerI_\h
\big( \hbox{\it co}\hskip1pt\Cal{S} \big) $  are again lattices w.r.t.~set
theoretical inclusion, so they can (and will) be thought of as categories
as well.

\vskip9pt

   {\bf Remarks 2.8:}  {\it (a)}  \, Let  $ \, X \in \{ \calI, \calC,
\gerI, \gerC \} \, $  and  $ \, S_\h \in \big\{ \fhg, \uhg \big\} \, $.
Since $ \, \pi_{S_\h}(X_\h) = X_\h \Big/ \big( X_\h \cap \h \, S_\h \big)
\, $, \, the property  $ \, X_\h \Big/ \h \, X_\h \cong \pi_{S_\h}(X_\h)
= X \, $  is equivalent to  $ \, X_\h \cap \h \, S_\h = \h \, X_\h \, $.
Therefore our quantum objects can also be characterized, instead of by
(2.3), by
 \vskip-9pt
  $$  \hbox{ $ \matrix
   (a) \qquad  &  \calI_\h \ideal_\ell \coideal \; \fhg \quad ,
&  \quad  \calI_\h \cap \h \, \fhg \, = \, \h \, \calI_\h \quad ,
&  \quad  \calI_\h \big/ \h \, \calI_\h \, = \, \calI  \\
   (b) \qquad  &  \calC_\h \leq^1 \coideal_\ell \, \fhg \quad ,
&  \quad  \calC_\h \cap \h \, \fhg \, = \, \h \, \calC_\h \quad ,
&  \quad  \calC_\h \big/ \h \, \calC_\h \, = \, \calC  \\
   (c) \qquad  &  \gerI_\h \ideal_\ell \coideal \;\; U_\h(\gerg) \quad ,
&  \quad  \gerI_\h \cap \h \, U_\h(\gerg) \, = \, \h \, \gerI_\h \quad ,
&  \quad  \gerI_\h \big/ \h \, \gerI_\h \, = \, \gerI  \\
   (d) \qquad  &  \gerC_\h \leq^1 \coideal_\ell \, U_\h(\gerg) \quad ,
&  \quad  \gerC_\h \cap \h \, U_\h(\gerg) \, = \, \h \, \gerC_\h \quad ,
&  \quad  \gerC_\h \big/ \h \, \gerC_\h \, = \, \gerC
               \endmatrix $ }   \eqno (2.3)'  $$
along with conditions (2.4).  In any case, next Lemma proves that the
formal subgroup of  $ G $  obtained as specialization of a quantum
formal subgroup is always coisotropic (much like specializing a
quantum group one gets a  {\sl Poisson\/}  group).
                                                  \par
   {\it (b)}  \, If a quadruple  $ \, \big( \calI_\h \, , \, \calC_\h
\, , \, \gerI_\h \, , \, \gerC_\h \big) \, $  is given which enjoys
all properties in the first and the second column of  $ (2.3)' $,
then one easily checks that the four specialized objects  $ \, \calI
:= \calI_\h\big|_{\h=0} \, $, $ \, \calC := \calC_\h\big|_{\h=0}
\, $,  $ \, \gerI := \gerI_\h\big|_{\h=0} \, $  and  $ \, \gerC :=
\gerC_\h\big|_{\h=0} \, $  verify relations  {\it (1)\/}  and  {\it
(2)\/}  in \S 1.2, thus they define one single pair  {\sl (coisotropic
subgroup, Poisson quotient)},  and the quadruple  $ \, \big( \calI_\h \, ,
\, \calC_\h \, , \, \gerI_\h \, , \, \gerC_\h \big) \, $  then yields a
quantization of the latter in the sense of \S 2.7.
                                                  \par
   {\it (c)}  \, The  {\sl existence\/}  of quantizations for a given
formal coisotropic subgroup is an open question, in general.  However,
Etingof and Kahzdan provided a positive answer for the special subclass
of those formal coisotropic subgroups  $ K $  which are also  {\sl Poisson
subgroups\/}  (which infinitesimally amounts to  $ \, \gerk := \text{\it
Lie}(K) \, $  being a Lie subbialgebra); see [EK2, \S 2.2].  Several
other examples of quantizations exist in literature for scattered cases
of special coisotropic subgroups of interest: we shall deal with one of
them in \S 6.   $ \quad \diamondsuit $

\vskip7pt

\proclaim{Lemma 2.9} Let  $ K $  be a formal subgroup of  $ G $,  and
assume a quantization  $ \, \calI_\h $,  $ \calC_\h $,  $ \gerI_\h $
or  $ \, \gerC_\h $  of  $ \, \calI $,  $ \calC $,  $ \gerI $  or
$ \, \gerC $  respectively be given as in \S 2.7.  Then  $ K $
is coisotropic.
\endproclaim

\demo{Proof} Assume  $ \, \calI_\h $  exists.  Let  $ \, f, g \in \calI
\, $,  \, and let  $ \, \varphi, \gamma \in \calI_\h \, $  with  $ \,
\pi_{F_\h}(\varphi) = f \, $,  $ \, \pi_{F_\h}(\gamma) = g \, $.  Then
by definition  $ \, \{f,g\} = \pi_{F_\h}\big( \h^{-1} [\varphi,
\gamma] \big) \, $.  But  $ \, [\varphi,\gamma] \in \calI_\h \cap \h
\, F_\h[G] = \h \, \calI_\h \, $  by assumption, hence  $ \, \h^{-1}
[\varphi,\gamma] \in \calI_\h \, $,  \, thus  $ \, \{f,g\} = \pi_{F_\h}
\big( \h^{-1} [\varphi,\gamma] \big) \in \pi_{F_\h}(\calI_\h) = \calI
\, $,  \, which means that  $ \, \calI \, $  is closed for the Poisson
bracket.  Thus (see \S 1.6)  $ K $  is coisotropic.  The proof is
entirely similar when dealing with  $ \calC_\h \, $,  $ \gerI_\h $  or
$ \gerC_\h \, $.   \qed
\enddemo

\vskip7pt

   {\bf \hskip-8pt  2.10 General program.} Starting from the 
       \hbox{setup of \S 1.2, we will move along the scheme}
 \vskip-7pt
  $$  \leqalignno{
   \qquad \;  \calI \; \big( \! \subseteq \! F[[G]] \big)  \;{\buildrel
(1) \over \longrightarrow}\;  &  \calI_\h \; \big( \! \subseteq \!
F_\h[[G]] \big)  \;{\buildrel (2) \over \longrightarrow}\;
{\calI_\h}^{\!\curlyvee} \, \big( \! \subseteq \! {F_\h[[G]]}^\vee \big)
\;{\buildrel (3) \over \longrightarrow}\;  {\calI_0}^{\!\curlyvee} \,
\Big( \! \subseteq \! {\big({F_\h[[G]]}^\vee\big)}_0 \! =
U\big(\gerg^*\big) \! \Big)  &  \hbox{\it (a)}  \cr
   \qquad \;  \calC \; \big( \! \subseteq \! F[[G]] \big)  \;{\buildrel
(1) \over \longrightarrow}\;  &  \calC_\h \; \big( \! \subseteq \!
F_\h[[G]] \big) \;{\buildrel (2) \over \longrightarrow}\;  {\calC_\h}^{\!\!
\triangledown} \, \big( \! \subseteq \! {F_\h[[G]]}^\vee \big)
\;{\buildrel (3) \over \longrightarrow}\;  {\calC_0}^{\!\!\triangledown}
\, \Big( \! \subseteq \! {\big({F_\h[[G]]}^\vee\big)}_0 \! =
U\big(\gerg^*\big) \! \Big)  &  \hbox{\it (b)}  \cr
   \qquad \,  \gerI \;\, \big( \! \subseteq U(\gerg) \big)  \;{\buildrel (1)
\over \longrightarrow}\;  &  \, \gerI_\h \; \big( \! \subseteq U_\h(\gerg)
\big)  \;{\buildrel (2) \over \longrightarrow}\; \, {\gerI_\h}^{\! !} \;
\big( \! \subseteq {U_\h(\gerg)}' \,\big) \;{\buildrel (3) \over
\longrightarrow}\;\,  {\gerI_0}^{\! !} \; \Big( \! \subseteq
{\big({U_\h(\gerg)}'\big)}_0 = F[[G^*]] \Big)
&  \hbox{\it (c)}  \cr
   \qquad \,  \gerC \;\, \big( \! \subseteq U(\gerg) \big)  \;{\buildrel
(1) \over \longrightarrow}\;  &  \, \gerC_\h \; \big( \! \subseteq
U_\h(\gerg) \big)  \;{\buildrel (2) \over \longrightarrow}\;\,
{\gerC_\h}^{\!\!\Lsh} \; \big( \! \subseteq {U_\h(\gerg)}' \,\big)
\;{\buildrel (3) \over \longrightarrow}\;\,  {\gerC_0}^{\!\!\Lsh} \;
\Big( \! \subseteq {\big({U_\h(\gerg)}'\big)}_0 = F[[G^*]] \Big)
&  \hbox{\it (d)}  \cr }  $$

\vskip2pt

   In the frame above, the arrows (1) are quantizations, as in \S 2.7,
and the arrows (3) are specializations at  $ \, \h = 0 \, $.  The middle
arrows (2) instead are suitable ``adaptations'' of Drinfeld's functors
to the quantizations of  $ K $  or of  $ G \big/ K $  in left hand side:
roughly, one takes the suitable Drinfeld's functor on  $ F[[G]] \, $,
resp.~on  $ U(\gerg) \, $,  and restricts it   --- in some sense ---
to the subobject  $ \calI $  or  $ \calC \, $,  resp.~$ \gerI $  or
$ \gerC \, $.  The points to show then are the following:
 \vskip3pt
   {\sl  $ \underline{\text{First}} $:} \, each one of the
right-hand-side objects above is one of the four algebraic
objects which describe a (closed formal) subgroup of  $ G^* \, $:
namely, the correspondence is
 \vskip1pt
   \centerline{ $
\text{\it (a)} \, =\joinrel\Longrightarrow \, \text{\it (c)}
\; ,  \hskip17pt
\text{\it (b)} \, =\joinrel\Longrightarrow \, \text{\it (d)}
\; ,  \hskip17pt
\text{\it (c)} \, =\joinrel\Longrightarrow \, \text{\it (a)}
\; ,  \hskip17pt
\text{\it (d)} \, =\joinrel\Longrightarrow \, \text{\it (b)}
\; . $ }
 \vskip1pt
   {\sl  $ \underline{\text{Second}} $:} \, all the formal subgroups of
$ G^* $  associated to the four objects so obtained are  {\sl coisotropic}.
                                       \par
   {\sl  $ \underline{\text{Third}} $:} \, the four formal subgroups of
$ G^* $  in  {\it (b)\/}  do coincide.
                                       \par
   {\sl  $ \underline{\text{Fourth}} $:} \, if we start from  $ \, K
\in \hbox{\it co}\hskip1pt\Cal{S}(G) \, $,  \, then the formal coisotropic
subgroup of  $ \, G^* $  obtained above is  $ K^\perp $  (cf.~Definition
1.4{\it (a)\/}).

\vskip1,3truecm

\centerline {\bf \S\; 3 \ Drinfeld-like functors on quantum
subgroups and Poisson quotients }

\vskip11pt

   In this section and next one we introduce Drinfeld-like functors for quantum coisotropic subgroups and Poisson quotients.  In particular, we start with
$ \, \calI_\h \, $,  $ \, \calC_\h \, $,  $ \, \gerI_\h \, $  and  $ \,
\gerC_\h \, $  as in \S 2.7, hence enjoying (2.3), or equivalently  $ (2.3)' $, and (2.4), with  $ F_\h $  and  $ U_\h $  as in \S 2.7.  We begin moving step (2) in \S 2.10, with a definition whose meaning is (roughly) to ``restrict'' Drinfeld's functors from quantum groups to quantum subgroups or Poisson quotients:

\vskip7pt

\proclaim {Definition 3.1}  {\sl (Drinfeld-like functors for subgroups)}
Keeping notation of \S 2.4, we define:
 \vskip-2pt
   \phantom{A}  (a) \hskip7pt  $ \displaystyle{ {\calI_\h}^{\!\curlyvee}
\; := \;  {\textstyle \sum_{n=1}^\infty} \, \h^{-n} \cdot I^{n-1}
\cdot \calI_\h  \; = \;  {\textstyle \sum_{n=1}^\infty} \, \h^{-n}
\cdot J^{n-1} \cdot \calI_\h  \quad ; } $
 \vskip1pt
   \phantom{A}  (b) \hskip7pt  $ \displaystyle{ {\calC_\h}^{\!\!
\triangledown}  \; := \;  \calC_\h \, + \, {\textstyle \sum_{n=1}^\infty}
\, \h^{-n} \cdot {\big( \calC_\h \cap I \,\big)}^n \, =  \;  \kh \cdot 1
\, + \, {\textstyle \sum_{n=1}^\infty} \, \h^{-n} \cdot {\big( \calC_\h
\cap J \,\big)}^n  \quad ; } $
 \vskip1pt
   \phantom{A}  (c) \hskip7pt  $ \displaystyle{ {\gerI_\h}^{\! !}  \;
:= \;  \Big\{\, x \in \gerI_\h \,\Big\vert\, \delta_n(x) \in \h^n \,
{\textstyle \sum_{s=1}^n} \, {U_\h}^{\otimeshat (s-1)} \otimeshat
\gerI_\h \otimeshat {U_\h}^{\otimeshat (n-s)} \! , \; \forall\; n
\in \N_+ \Big\}  \quad ; } $
 \vskip1pt
   \phantom{A}  (d) \hskip7pt  $ \displaystyle{ {\gerC_\h}^{\!\!\Lsh}
\; := \;  \Big\{\, x \in \gerC_\h \;\Big\vert\; \delta_n(x) \in \h^n
\, {U_\h}^{\otimeshat (n-1)} \otimeshat \gerC_\h \, , \; \forall\;
n \in \N_+ \,\Big\}  \quad . } $
\endproclaim

\vskip7pt

   {\bf 3.2 Remark:} \, The following inclusion relations hold, directly
by definitions:
  $$  (i) \quad\;  {\calI_\h}^{\!\curlyvee} \supseteq \calI_\h \; ,  \qquad
(ii) \quad\; {\calC_\h}^{\!\!\triangledown} \supseteq \calC_\h \; ,
\qquad  (iii) \quad\;  {\gerI_\h}^{\! !} \subseteq \gerI_\h \; ,
\qquad  (iv) \quad\;  {\gerC_\h}^{\!\!\Lsh} \subseteq \gerC_\h \; .  $$
 \eject
   Moreover, definitions and assumptions in  $ (2.3)' $  imply that
$ \, \calI_\h = {\calI_\h}^{\!\curlyvee} \cap F_\h \, $,  $ \, \calC_\h
= {\calC_\h}^{\!\!\triangledown} \cap F_\h \, $,  $ \, {\gerI_\h}^{\! !}
= \gerI_\h \cap {U_\h}' \, $  and  $ \, {\gerC_\h}^{\!\!\Lsh} = \gerC_\h
\cap {U_\h}' \, $:  \, thus we are just ``restricting'' Drinfeld's
functors.   $ \quad \diamondsuit $

\vskip7pt

   We can now state the QDP for formal coisotropic subgroups and
Poisson quotients:

\vskip7pt

\proclaim {Theorem 3.3} ({\sl ``QDP for Coisotropic Subgroups and
Poisson Quotients''})
                                        \par
   (a) \, Definition 3.1 provides category equivalences
 \vskip-15pt
  $$  \matrix
   {(\ )}^\curlyvee \colon \, \calI_\h \big( \hbox{\it co}
\hskip1pt\Cal{S}(G) \big) \, {\buildrel \cong \over \llongrightarrow}
\, \gerI_\h \big( \hbox{\it co}\hskip1pt\Cal{S}(G^*) \big) \; ,
&  \quad {(\ )}^\triangledown \colon \, \calC_\h \big( \hbox{\it
co}\hskip1pt\Cal{S}(G) \big) \, {\buildrel \cong \over \llongrightarrow}
\, \gerC_\h \big( \hbox{\it co}\hskip1pt\Cal{S}(G^*) \big)
\; ,  \\
   {(\ )}^! \, \colon \, \gerI_\h \big( \hbox{\it co}\hskip1pt
\Cal{S}(G) \big) \, {\buildrel \cong \over \llongrightarrow} \,
\calI_\h \big( \hbox{\it co}\hskip1pt\Cal{S}(G^*) \big) \; ,
&  \quad  {(\ )}^\Lsh \, \colon \, \gerC_\h \big( \hbox{\it co}
\hskip1pt\Cal{S}(G) \big) \, {\buildrel \cong \over \llongrightarrow}
\, \calC_\h \big( \hbox{\it co}\hskip1pt \Cal{S}(G^*) \big) \; ,
\endmatrix  $$
 \vskip-7pt
\noindent
 along with the similar ones with  $ G $  and  $ G^* $  interchanged,
such that  $ \, {(\ )}^! \circ {(\ )}^\curlyvee = \hbox{\sl id}_{co
\hskip1pt\Cal{S}(G^)} \, $,  $ \, {(\ )}^\curlyvee \circ {(\ )}^! =
\hbox{\sl id}_{co\hskip1pt\Cal{S}(G^*)} \, $,  \, and  $ \, {(\ )}^\Lsh
\circ {(\ )}^\triangledown = \hbox{\sl id}_{co\hskip1pt\Cal{S}(G)} \, $,
$ \, {(\ )}^\triangledown \circ {(\ )}^\Lsh = \hbox{\sl id}_{co\hskip1pt
\Cal{S}(G^*)} \, $,  \, and so on.
 \vskip7pt
   (b) \,  {\sl (QDP)}  For any  $ \, K \in
\hbox{\it co}\hskip1pt\Cal{S}(G) \, $,  \, we have
 \vskip-3pt
  $$  \matrix
   \calI(\gerk)_\h^{\,\curlyvee} \mod \h \, {F_\h[[G]]}^\vee
\; = \; \gerI\big(\gerk^\perp\big) \, ,  &
\quad  \calC(\gerk)_\h^{\,\triangledown} \mod \h \, {F_\h[[G]]}^\vee
\; = \; \gerC\big(\gerk^\perp\big) \, ,  \\
   \gerI(\gerk)_\h^{\;!} \mod \h \, {U_\h(\gerg)}'
\; = \; \calI\big(\gerk^\perp\big) \, ,  &
\quad  \gerC(\gerk)_\h^{\,\Lsh} \mod \, \h \, {U_\h(\gerg)}'
\; = \; \calC\big(\gerk^\perp\big) \, .
      \endmatrix  $$
In short, the quadruple  $ \, \big( \calI(\gerk)_\h^{\,\curlyvee},
\, \calC(\gerk)_\h^{\,\triangledown}, \, \gerI(\gerk)_\h^{\;!} \, ,
\, \gerC(\gerk)_\h^{\,\Lsh} \big) \, $  is a quantization of the
quadruple  $ \Big( \gerI\big(\gerk^\perp\big), \gerC\big(\gerk^\perp\big),
\calI\big(\gerk^\perp\big), \calC\big(\gerk^\perp\big) \!\Big) $
w.r.t.~the quantization  $ \big( {F_\h[[G]]}^\vee \! , {U_\h(\gerg)}'
\big) $  of  $ \, \big( U(\gerg^*), F[[G^*]] \big) $.
\endproclaim

\vskip1,3truecm

\centerline {\bf \S\; 4 \ First properties of Drinfeld-like functors }

\vskip11pt

   We shall now study the properties of the images of Drinfeld-like functors
for general  $ \h \, $.  The main result is   --- Proposition 4.4 ---   that
they are quantizations of some (unique) pair  {\sl (coisotropic subgroup,
Poisson quotient)},
   \hbox{in the sense of \S 2.7, for the Poisson group  $ G^* $.}

\vskip7pt

\proclaim {Lemma 4.1}  The following relations hold (w.r.t.~the perfect Hopf pairing between  $ {U_\h}' $  and  $ {F_\h}^{\!\vee} $  given by Proposition 2.6 for the orthogonality relations  {\it (i)--(ii)\/}):
 \vskip-13pt
  $$  \matrix
   \hbox{\it (i)}  \;\quad  {\calI_\h}^{\!\curlyvee} =
{\big( {\gerC_\h}^{\!\!\Lsh\,} \big)}^\perp \, ,  \;\quad
{\gerC_\h}^{\!\!\Lsh} = {\big( {\calI_\h}^{\!\curlyvee} \big)}^\perp
\hskip7pt
  &  \hbox{\it (ii)}  \;\quad  {\gerI_\h}^{\! !} =
{\big( {\calC_\h}^{\!\!\triangledown} \big)}^\perp \, ,  \;\quad
{\calC_\h}^{\!\!\triangledown} = {\big( {\gerI_\h}^{\! ! \,}
\big)}^\perp  \\
   \hbox{\it (iii)}  \quad  {\calI_\h}^{\!\curlyvee} =
{F_\h}^{\!\vee} \cdot {\big( {\calC_\h}^{\!\!\triangledown} \big)}^{\!+}
\, ,  \;\;  {\calC_\h}^{\!\!\triangledown} = {\big( {F_\h}^{\!\vee}
\,\big)}^{\text{\it co}{\calI_\h}^{\!\curlyvee}}  \hskip1pt
  &  \hbox{\it (iv)}  \quad  {\gerI_\h}^{\! !} =
{U_\h}^{\!\prime} \cdot {\big( {\gerC_\h}^{\!\!\Lsh} \,\big)}^{\!+}
\, ,  \;\;  {\gerC_\h}^{\!\!\Lsh} = {\big( {U_\h}^{\!\prime}
\,\big)}^{\text{\it co}{\gerI_\h}^{\! !}}  \\
   \endmatrix  $$
\endproclaim

\demo{Proof}  Let  $ \, I = I_{\scriptscriptstyle F_\h} \, $  be
the ideal of  $ F_\h $  considered in \S 2.4, and take  $ \, y_1 $,
$ \dots $,  $ y_{n-1} \in I \, $;  then  $ \, \langle y_i, 1 \rangle
= \epsilon(y_i) \in \h \cdot \kh \, $,  \, for all  $ \, i= 1, \dots,
n-1 \, $.  Given  $ \, y_n \in \calI_\h \, $  and  $ \, \gamma \in
{\gerC_\h}^{\!\!\Lsh} \, $,  consider
  $$  \left\langle \, {\textstyle \prod\limits_{i=1}^n} y_i \, , \, \gamma
\right\rangle = \left\langle \, {\textstyle \bigotimes\limits_{i=1}^n}
y_i \, , \Delta^n(\gamma) \right\rangle = \bigg\langle \, {\textstyle
\bigotimes\limits_{i=1}^n} y_i \, , {\textstyle \sum\limits_{\Psi
\subseteq \{1,\dots,n\}}} \! \delta_\Psi(\gamma) \bigg\rangle = \!
{\textstyle \sum\limits_{\Psi \subseteq \{1,\dots,n\}}} \!
\bigg\langle \, {\textstyle \bigotimes\limits_{i=1}^n} c_i \, ,
\delta_\Psi(\gamma) \bigg\rangle \, .  $$
   \indent   Now consider any summand in the last term in the formula
above.  Let  $ \, \vert\Psi\vert = t \, $  ($ \, t \leq n \, $):  then
$ \, \big\langle \otimes_{i=1}^n y_i, \delta_\Psi(\gamma) \big\rangle
= \big\langle \otimes_{i \in \Psi} y_i, \delta_t(\gamma) \big\rangle
\cdot \prod_{j \not\in \Psi} \langle y_j, 1 \rangle \, $,  by
definition of  $ \, \delta_\Psi \, $.  Thanks to the previous
analysis, we have  $ \, \prod_{j \not\in \Psi} \langle y_j, 1
\rangle \in \h^{n-t} \kh \, $,  \, hence
  $$  \Big\langle \, {\textstyle \bigotimes_{i \in \Psi}} y_i \, ,
\, \delta_t(\gamma) \Big\rangle \in \left\langle \, {\textstyle
\bigotimes_{i \in \Psi}} y_i \, , \, \h^t \, {\textstyle \sum_{s=1}^n}
\, {U_\h}^{\otimeshat (n-1)} \otimeshat \gerC_\h \right\rangle \subseteq
\h^{t+1} \, \kh  $$
because  $ \, \gamma \in {\gerC_\h}^{\!\!\Lsh} \, $;  therefore  $ \,
\left\langle \, {\textstyle \prod_{i=1}^n} y_i \, , \, \gamma \right\rangle
\in \h \, \kh \, $.  And even more, the rightmost tensor factor in each
summand  $ \delta_\Psi(\gamma) $  always belongs to  $ \gerC_\h $  (as
also  $ \, 1 \in \gerC_\h \, $),  whereas  $ \, y_n \in \calI_\h =
{\gerC_\h \phantom{)}}^{\!\!\!\!\perp} \, $:  \, therefore  $ \,
\big\langle \prod_{i=1}^n y_i \, , \, \gamma \big\rangle = \Big\langle
\, {\textstyle \bigotimes_{i=1}^n} y_i \, , \sum_{\Psi \subseteq \{1,
\dots,n\}} \! \delta_\Psi(\gamma) \Big\rangle = 0 \, $.  This means that
  $$  {\calI_\h}^{\!\curlyvee} \subseteq {\big( {\gerC_\h}^{\!\!\Lsh\,}
\big)}^\perp \; ,  \qquad \qquad  {\gerC_\h}^{\!\!\Lsh} \subseteq
{\big( {\calI_\h}^{\!\curlyvee} \big)}^\perp \; .   \eqno (4.1)  $$
   \indent   Now take  $ \, \kappa \in {\big( {\calI_\h}^{\!\curlyvee}
\big)}^\perp \subseteq {\big({F_\h}^{\!\vee}\big)}^* = {U_\h}' \, $
(using Proposition 2.6 for the last equality).  Since  $ \, \kappa
\in {U_\h}' \, $,  we have  $ \, \delta_n(\kappa) \in \h^n \,
{U_\h}^{\otimeshat n} \, $  for all  $ \, n \in \N \, $,  and
moreover from  $ \, \kappa \in {\big( {\calI_\h}^{\!\curlyvee}
\big)}^\perp \, $  it follows that  $ \, \kappa_+ := \h^{-n} \,
\delta_n(\kappa) \, $  enjoys  $ \, \Big\langle I^{\otimestilde
(n-1)} \otimestilde \calI_\h \, , \, \kappa_+ \Big\rangle = 0 \, $,
\, so that
  $$  \kappa_+ \in {\Big( I^{\otimestilde (n-1)} \otimestilde \calI_\h
\Big)}^{\!\perp} = \; {\textstyle \sum_{r+s=n-2}} \, {U_\h}^{\otimeshat r}
\otimeshat I^\perp \otimeshat {U_\h}^{\otimeshat s} \otimeshat U_\h \,
+ \, {U_\h}^{\otimeshat (n-1)} \otimeshat {\calI_\h \phantom{)}}^{\!\!
\!\perp} \;\; .  $$
In addition,  $ \, \delta_n(\kappa) \in J^{\otimeshat n} \, $,  where
$ \, J := J_{U_\h} = \Ker\,\big( \, \epsilon \! : U_\h \longrightarrow
\kh \big) \, $,  hence  $ \, \delta_n(\kappa) \in \h^n \,
{U_\h}^{\otimeshat n} \cap J^{\otimeshat n} = \h^n \, J^{\otimeshat n}
\, $;  \, this together with the above formula yields
  $$  \displaylines{
   \kappa_+ \in {\Big( I^{\otimestilde (n-1)} \otimestilde \calI_\h
\Big)}^{\!\perp} \cap J^{\otimeshat n}  =  \bigg( {\textstyle
\sum\limits_{r+s=n-2}} {U_\h}^{\widehat\otimes\, r} \otimes I^\perp
\otimeshat {U_\h}^{\widehat\otimes\, s} \otimeshat U_\h \bigg) \cap
J^{\widehat\otimes\, n} \, +   \hfill  \cr
   \hfill   + \, \Big( {U_\h}^{\widehat\otimes
\, (n-1)} \otimeshat {\calI_\h \phantom{(}}^{\!\!\!\!\perp} \Big)
\cap J^{\widehat\otimes\, n}  =  \hskip-3pt  {\textstyle
\sum\limits_{r+s=n-2}}  \hskip-3pt  J^{\widehat\otimes\, r}
\otimeshat \Big( I^\perp \cap J_U \Big) \! \otimeshat
J^{\widehat\otimes\, s} \otimeshat J \, + \,
J^{\widehat\otimes\, (n-1)} \otimeshat \Big( \calI_\h
\phantom{(}^{\!\!\!\!\perp} \cap J \Big) =  \cr
   \hfill   = J^{\widehat\otimes\, (n-1)} \otimeshat \Big( {\calI_\h
\phantom{(}}^{\!\!\!\perp} \cap J \Big) = J^{\widehat\otimes\, (n-1)}
\otimeshat \Big( \gerC_\h \cap J \Big) \subseteq {U_\h}^{\widehat\otimes\,
(n-1)} \otimeshat \gerC_\h  \cr }  $$
where in the third equality we used the fact that  $ \, I^\perp = 0 \, $;
\, the last equality then follows from  (2.4){\it (i)}.  Thus  $ \,
\kappa_+ \in {U_\h}^{\widehat\otimes\, (n-1)} \otimeshat \gerC_\h \, $,
hence  $ \, \delta_n(\kappa) \in \h^n \, {U_\h}^{\widehat\otimes\, (n-1)}
\otimeshat \gerC_\h \, $  for all  $ \, n \in \N \, $:  so  $ \, \kappa
\in {\gerC_\h}^{\!\!\Lsh} \, $.  We conclude that  $ \, {\big(
{\calI_\h}^{\!\curlyvee} \big)}^\perp \subseteq {\gerC_\h}^{\!\!\Lsh}
\, $, which together with  $ (4.1) $  gives  $ \; {\gerC_\h}^{\!\!\Lsh}
= {\big( {\calI_\h}^{\!\curlyvee} \big)}^\perp \, $.
                                        \par
   By Proposition 2.6 the specialization at  $ \, \h = 0 \, $  of the
pairing between  $ {U_\h}' $  and  $ {F_\h}^{\!\vee} $  is perfect too.
From this we can easily argue that  $ \, {\calI_\h}^{\!\curlyvee} \equiv
\Big(\! {\big( {\calI_\h}^{\!\curlyvee} \big)}^\perp \Big)^{\!\perp} \!
\mod \h \, {F_\h}^{\!\vee} \, $,  \, whence  $ \, {\calI_\h}^{\!\curlyvee}
= \Big(\! {\big( {\calI_\h}^{\!\curlyvee} \big)}^\perp \Big)^{\!\perp}
\, $  follows at once by  $ \h $--adic  completeness.  But then
starting from  $ \, {\gerC_\h}^{\!\!\Lsh} = {\big( {\calI_\h}^{\!
\curlyvee} \big)}^\perp \, $,  \, hence  $ \, {\big( {\gerC_\h}^{\!
\!\Lsh\,} \big)}^\perp = \Big(\! {\big( {\calI_\h}^{\!\curlyvee}
\big)}^\perp \Big)^{\!\perp} \, $,  \, we finally get  $ \, {\big(
{\gerC_\h}^{\!\!\Lsh\,} \big)}^\perp = {\calI_\h}^{\!\curlyvee} \, $,
\, thus  {\it (i)\/}  is proved.
                                        \par
   The proof of  {\it (ii)}  is similar.  First of all, by (2.4){\it
(ii)\/}  and definitions it is clear that
  $$  {\gerI_\h}^{\! !} \subseteq {\big( {\calC_\h}^{\!\!\triangledown}
\big)}^\perp \, ,  \qquad \qquad  {\calC_\h}^{\!\!\triangledown} \subseteq
{\big( {\gerI_\h}^{\! ! \,} \big)}^\perp \, .   \eqno (4.2)  $$
   \indent   Now notice that  $ \, {\calC_\h}^{\!\!\triangledown} \supseteq
\calC_\h \, $,  so  $ \, {\big( {\calC_\h}^{\!\!\triangledown} \big)}^\perp
\subseteq {\calC_\h \phantom{)}}^{\!\!\!\!\perp} = \gerI_\h \, $,  due to
(2.4){\it (ii)\/};  thus  $ \, {\big( {\calC_\h}^{\!\!\triangledown}
\big)}^\perp \subseteq \gerI_\h \, $.   Second, pick  $ \, \eta \in
{\big( {\calC_\h}^{\!\!\triangledown} \big)}^\perp \; \big(
\subseteq {U_\h}' \big) \, $.  Then  $ \, \delta_n(\eta) \in \h^n
{U_\h}^{\widehat\otimes\, n} \, $  for all  $ \, n \in \N_+ \, $,
\, and from  $ \, \eta \in {\big( {\calC_\h}^{\!\!\triangledown}
\big)}^\perp \, $  we get that  $ \, \eta_+ := \h^{-n} \,
\delta_n(\eta) \, $  enjoys  $ \, \Big\langle {\big( \calC_\h
\cap I \,\big)}^{\widetilde\otimes\, n} , \, \eta_+ \Big\rangle
= 0 \; $,  \; so that
  $$  \eta_+ \, \in \, \Big(\! {\big( \calC_\h \cap I
\,\big)}^{\widetilde\otimes \, n} \Big)^{\!\perp} = \,
{\textstyle \sum\limits_{r+s=n-1}} {U_\h}^{\widehat\otimes\, r}
\otimeshat {\big( \calC_\h \cap I \,\big)}^{\!\perp} \otimeshat
{U_\h}^{\widehat\otimes\, s} \; .  $$
   \indent   Moreover  $ \, \delta_n(\eta) \in J^{\widehat\otimes\, n}
\, $,  hence  $ \, \delta_n(\eta) \in \h^n \, {U_\h}^{\widehat\otimes
\, n} \cap \, J^{\widehat\otimes\, n} = \, \h^n J^{\widehat\otimes\, n}
\, $,  \; so  $ \, \eta_+ \in J^{\widehat\otimes\, n} \, $  and
  $$  \displaylines{
   \eta_+ \; \in \; \Big(\! {\big( \calC_\h \cap I
\,\big)}^{\widehat\otimes\, n} \Big)^{\!\perp} \cap
J^{\widehat\otimes\, n} \, = \, \Big( {\textstyle \sum_{r+s=n-1}}
\, {U_\h}^{\widehat\otimes\, r} \otimeshat {\big( \calC_\h \cap I
\,\big)}^{\!\perp} \otimeshat {U_\h}^{\widehat\otimes\, s} \Big)
\cap J^{\widehat\otimes\, n} \, =  \cr
   \hfill   = \, {\textstyle \sum_{r+s=n-1}} \, J^{\widehat\otimes\, r}
\otimeshat \Big( \! {\big( \calC_\h \cap I \,\big)}^{\!\perp} \cap J \Big)
\otimeshat J^{\widehat\otimes\, s} \; .  \cr }  $$
Now  $ \, {\big( \calC_\h \cap I \,\big)}^\perp \cap J = {\calC_\h
\phantom{)}}^{\!\!\!\perp} \cap J = \gerI_\h \cap J \subseteq
\gerI_\h \, $,  \; thanks to  (2.4){\it (ii)}.  The upshot is
  $$  \eta_+ \; \in \; {\textstyle \sum_{r+s=n-1}} \, J^{\widehat\otimes
\, r} \otimeshat \big( \gerI_\h \cap J_U \,\big) \otimeshat
J^{\widehat\otimes\, s} \, \subseteq \, {\textstyle \sum_{r+s=n-1}}
\, {U_\h}^{\widehat\otimes\, r} \otimeshat \gerI_\h \otimeshat
{U_\h}^{\widehat\otimes\, s}  $$
whence we get  $ \, \delta_n(\eta) \in \h^n \sum_{r+s=n-1}
\, {U_\h}^{\widehat\otimes\, r} \otimeshat \gerI_\h \otimeshat
{U_\h}^{\otimeshat s} \, $  for all  $ \, n \in \N_+ \, $.  Since
in addition  $ \, \eta \in \gerI_\h \, $,  for we proved that  $ \,
{\big( {\calC_\h}^{\!\!\triangledown} \big)}^\perp \subseteq \gerI_\h
\, $,  we argue that  $ \, \eta \in {\gerI_\h}^{\! !} \, $.  The final
outcome is  $ \, {\big( {\calC_\h}^{\!\!\triangledown} \big)}^\perp
\subseteq {\gerI_\h}^{\! !} \, $,  \, which together with (4.2) implies
$ \, {\gerI_\h}^{\! !} = {\big( {\calC_\h}^{\!\!\triangledown}
\big)}^\perp \, $,  \, q.e.d.
                                        \par
   With like arguments as for part  {\it (i)\/}  one proves that  $ \,
\Big(\! {\big( {\calC_\h}^{\!\!\triangledown} \big)}^\perp \Big)^{\!\perp}
\! = \, {\calC_\h}^{\!\!\triangledown} \, $  and then argue that  $ \,
{\big( {\gerI_\h}^{\! !} \,\big)}^\perp \! = \, {\calC_\h}^{\!\!
\triangledown} \, $;  \, this ends the proof of claim  {\it (ii)\/}
too.  Finally,  {\it (iii)\/}  and  {\it (iv)\/}  are straightforward
consequence of relations  {\it (iii)\/}  and  {\it (iv)\/}  in (2.4)
and of definitions.   \qed
\enddemo

\vskip7pt

\proclaim{Lemma 4.2}
  $$  \matrix
   \text{(a)} \;\; {\calI_\h}^{\!\curlyvee} \ideal_\ell {F_\h}^{\!\vee}
&  {} \quad\!  \text{(b)} \;\; {\calC_\h}^{\!\!\triangledown} \leq^1
\! {F_\h}^{\!\vee}  &  {} \quad\!  \text{(c)} \;\; {\gerI_\h}^{\! !}
\ideal_\ell {U_\h}'  &  {} \quad\!  \text{(d)} \;\;
{\gerC_\h}^{\!\!\Lsh} \leq^1 \! {U_\h}' \, ;  \\
   {}  \\
   \text{(e)} \;\; {\calI_\h}^{\!\curlyvee} \,\coideal\; {F_\h}^{\!\vee}
&  {} \quad\!  \text{(f)} \;\; {\calC_\h}^{\!\!\triangledown}
\,\coideal_\ell \, {F_\h}^{\!\vee}  &  {} \quad\!  \text{(g)} \;\;
{\gerI_\h}^{\! !} \,\coideal\; {U_\h}'  &  {} \quad\! \text{(h)}
\;\; {\gerC_\h}^{\!\!\Lsh} \,\coideal_\ell\, {U_\h}'  \\
      \endmatrix  $$
\endproclaim

\demo{Proof}  The statements on the first line are proved directly, and
imply those on the second line via the orthogonality relations of
Lemma 4.1.
                                          \par
   Claim  {\it (a)\/}  is straightforward, and  {\it (b)\/}  follows directly from definitions.  To prove  {\it (c)},  let  $ \, a \in {U_\h}' \, $  and
$ \, b \in {\gerI_\h}^{\! !} \, $:  by definition of  $ \, {\gerI_\h}^{\! !}
\, $,  from  $ \, \gerI_\h \ideal_\ell U_\h \, $  and from (2.1) we get  $ \,
\delta_n(a b) \in \h^n \sum_{s=1}^n {U_\h}^{\widehat\otimes\, (s-1)} \!
\otimeshat \gerI_\h \otimeshat {U_\h}^{\widehat\otimes\, (n-s)} $,  so
$ \, a b \in {\gerI_\h}^{\! !} \, $;  \, thus  $ \, {\gerI_\h}^{\! !}
\ideal_\ell {U_\h}' \, $.  Recall that  $ \, {U_\h}' \, $  is commutative
modulo  $ \h \, $,  and  $ \, \h \, {U_\h}' \in {\gerI_\h}^{\! !} \, $:
\, then  $ \, {\gerI_\h}^{\! !} \ideal_\ell {U_\h}' \, $  implies  $ \,
{\gerI_\h}^{\! !} \ideal {U_\h}' \, $  (a two-sided ideal), thus proving
{\it (c)}.  Lastly, to prove  {\it (d)},  remark that  $ \, 1 \in \gerC_\h \, $  and  $ \, \delta_n(1) = 0 \, $  for all  $ \, n \in \N \, $,  so  $ \, 1 \in {\gerC_\h}^{\!\!\Lsh} \, $.  Let  $ \, x, y \in {\gerC_\h}^{\!\!\Lsh} \, $  and  $ \, n \in \N \, $;  by (2.1) we have  $ \, \delta_n(x y) = \sum_{\Lambda \cup
Y = \{1,\dots,n\}}  \delta_\Lambda(x) \, \delta_Y(y) \, $.  Each of the factors
$ \, \delta_\Lambda(x) \, $  belongs to a module  $ \, \h^{\vert \Lambda \vert} \, {U_\h}^{\widehat\otimes\, ( \vert \Lambda \vert - 1 )} \! \otimeshat X \, $  where the last tensor factor is either  $ \, X = \gerC_\h \, $  (if  $ \, n \in \Lambda \, $)  or  $ \, X = \{1\} \subset \gerC_\h \, $  (if  $ \, n \not\in \Lambda \, $),  and similarly for  $ \, \delta_Y(y) \, $;  \, but
$ \, \Lambda \cup Y = \{1,\dots,n\} \, $  implies  $ \, \vert \Lambda \vert
+ \vert Y \vert \geq n \, $,  \, and summing up  $ \, \delta_n(x y) \in
\h^n {U_\h}^{\widehat\otimes\, (n-1)} \! \otimeshat \gerC_\h \, $,  whence
$ \, x y \in {\gerC_\h}^{\!\!\Lsh} \, $.  Thus  $ \, {\gerC_\h}^{\!\!\Lsh}
\leq^1 \! {U_\h}' \, $,  \, q.e.d.   \qed
\enddemo

 \vskip4pt

   {\it  $ \underline{\hbox{\it Remark}} $:} \, in the previous proof one
might also prove the required properties for only one of the objects involved, say  $ {\gerI_\h}^{\! !} $  for instance: then the properties of all others objects will follow from relations  {\it (i)--(iv)\/}  in Lemma 4.1.

 \vskip7pt

\proclaim{Lemma 4.3}
  $$  \matrix
   \hskip10pt  (a) \hskip7pt  {\calI_\h}^{\!\curlyvee} {\textstyle
\bigcap} \; \h \, {F_\h}^{\!\vee} \, = \; \h \, {\calI_\h}^{\!\curlyvee}
\, ,  &  \qquad  \hskip4pt  (b) \hskip7pt  {\calC_\h}^{\!\!\triangledown}
\, {\textstyle \bigcap} \; \h \, {F_\h}^{\!\vee} \, =
\; \h \, {\calC_\h}^{\!\!\triangledown} \, ,  \\
   (c) \hskip7pt  {\gerI_\h}^{\! !} \, {\textstyle \bigcap} \;
\h \, {U_\h}' \, = \; \h \, {\gerI_\h}^{\! !} \; ,  &  \qquad
   (d) \hskip7pt  {\gerC_\h}^{\!\!\Lsh} \, {\textstyle \bigcap} \;
\h \, {U_\h}' \, = \; \h \, {\gerC_\h}^{\!\!\Lsh} \; .
      \endmatrix  $$
\endproclaim
 \eject

\demo{Proof} We start proving claim  {\it (c)}.  Let  $ \; \eta \in
{\gerI_\h}^{\! !} \cap \h \, {U_\h}' \, = \, \h \, {\gerI_\h}^{\! !}
\, $.  Then
  $$  \delta_n(\eta) \, \in \, \h^n \left( \left( {\textstyle
\sum_{s=1}^n} {U_\h}^{\widehat\otimes\, (s-1)} \otimeshat \gerI_\h
\otimeshat {U_\h}^{\widehat\otimes\, (n-s)} \right) {\textstyle \bigcap}
\, \h \, {U_\h}^{\widehat\otimes\, n} \right)   \eqno (4.3)  $$
for all  $ \, n \in \N_+ \, $.  Now, for  $ \, n \in \N_+ \, $  we have
$ \; \left( {\textstyle \sum_{s=1}^n} {U_\h}^{\widehat\otimes\, (s-1)}
\otimeshat \gerI_\h \otimeshat {U_\h}^{\widehat\otimes\, (n-s)} \right)
{\textstyle \bigcap} \, \h \, {U_\h}^{\widehat\otimes\, n} = {\textstyle
\sum_{s=1}^n} \, {U_\h}^{\widehat\otimes \, (s-1)} \otimeshat \Big(
\gerI_\h \, {\textstyle \bigcap} \, \h \, U_\h \Big) \otimeshat
{U_\h}^{\widehat\otimes\, (n-s)} \, $,  \, and since  $ \, \gerI_\h
\, {\textstyle \bigcap} \, \h \, U_\h = \h \, \gerI_\h \, $  by
$ (2.3)' $,  from (4.3) we conclude that  $ \; \delta_n(\eta) \,
\in \, \h^{n+1} {\textstyle \sum_{s=1}^n} {U_\h}^{\widehat\otimes\,
(s-1)} \otimeshat \gerI_\h \otimeshat {U_\h}^{\widehat\otimes\, (n-s)}
\; $  for all  $ \, n \in \N_+ \, $,  \, which in turn means  $ \, \eta
\in \! \h \, {\gerI_\h}^{\! !} \, $,  \, q.e.d.  The converse inclusion
$ \, {\gerI_\h}^{\! !} \cap \h \, {U_\h}' \, \supseteq \, \h \,
{\gerI_\h}^{\! !} \, $  is trivially true.  The same arguments
prove  {\it (d)\/}  as well.
                                              \par
  As for  {\it (a)\/}  and  {\it (b)},  we can give a rather concrete
description of the objects involved, starting from  $ {F_\h}^{\!\vee} $.
Let  $ \, I := I_{\scriptscriptstyle F_\h} \, $  as in \S 2.4, $ \, J
:= \text{\it Ker}\, \big(\epsilon \, \colon \, F_\h \longrightarrow \kh
\,\big) \, $,  \, and  $ \, J^\vee := \h^{-1} J \subset {F_\h}^{\!\vee}
\, $.  Then  $ \, J \! \mod \h \, F_\h = J_{\scriptscriptstyle G} :=
\text{\it Ker}\, \big(\,\epsilon \, \colon \, F[[G]] \rightarrow
\Bbbk \,\big) \, $,  \, and  $ \, J_{\scriptscriptstyle G} \Big/
{J_{\scriptscriptstyle G}}^{\!2} = \gerg^* \, $.  Let  $ \, \{y_1,
\dots, y_n\} \, $,  with  $ \, n := \dim(G) \, $,  be a  $ \Bbbk $--basis
of  $ \, J_{\scriptscriptstyle G} \Big/ {J_{\scriptscriptstyle G}}^{\!2}
\, $,  \, and pull it back to a subset  $ \, \{j_1,\dots,j_n\} \, $  of
$ J \, $.  Then  $ \, \big\{ \h^{-|\underline{e}|} j^{\,\underline{e}}
\mod \h \, {F_\h}^{\!\vee} \;\big|\; \underline{e} \in \N^{\,n} \,\big\}
\, $  (with  $ \, j^{\,\underline{e}} := \prod_{s=1}^n j_s^{\,
\underline{e}(i)} \, $,  \, and similarly hereafter) is a
$ \Bbbk $--basis  of  $ {F_0}^{\!\vee} $  and, setting  $ \,
j_s^{\,\vee} := \h^{-1} j_s \, $  for all  $ s $,  the set
$ \, \big\{ j_1^{\,\vee}, \dots, j_n^{\,\vee} \big\} \, $
is a  $ \Bbbk $--basis  of  $ \; \gert := J^\vee \mod \h \,
{F_\h}^{\!\vee} \, $.  Moreover, since  $ \, j_\mu \, j_\nu
- j_\nu \, j_\mu \in \h \, J \, $  (for  $ \, \mu, \nu \in \{1,
\dots, n\} \, $)  we have  $ \; j_\mu \, j_\nu - j_\nu \, j_\mu
= \h \sum_{s=1}^n c_s \, j_s + \h^2 \gamma_1 + \h \, \gamma_2 \; $
for some  $ \, c_s \in \kh \, $,  $ \, \gamma_1 \in J \, $  and  $ \,
\gamma_2 \in J^2 $,  \, whence  $ \; \big[ j_\mu^\vee, j_\nu^\vee \,\big]
:= j_\mu^\vee \, j_\nu^\vee - j_\nu^\vee \, j_\mu^\vee
%
%
 \equiv {\textstyle \sum_{s=1}^n} \,
c_s \, j_s^\vee \; \mod \, \h \, {F_\h}^{\!\vee} \, $,
%
%
\, thus  $ \; \gert
:= J^\vee \! \mod \h \, {F_\h}^{\!\vee} \, $  is a Lie subalgebra of
$ {F_0}^{\!\vee} \, $:  \, indeed,  $ \, {F_0}^{\!\vee} = U(\gert) \, $
as Hopf algebras.
                                             \par
   Now for the second step.  The specialization map  $ \; \pi^\vee
\colon \, {F_\h}^{\!\vee} \relbar\joinrel\twoheadrightarrow
{F_0}^{\!\vee} = U(\gert) \; $  restricts to  $ \; \eta \, \colon
\, J^\vee \! \loongtwoheadrightarrow \gert \, := J^\vee
\!\!\! \mod \h \, {F_\h}^{\!\vee} = J^\vee \Big/ J^\vee \cap
{\big( \h \, {F_\h}^{\!\vee} \big)} = J^\vee \Big/ \big( J
+ J^\vee J_\h \big) \, $,  \, because  $ \, J^\vee \cap
{\big( \h \, {F_\h}^{\!\vee} \big)} = J^\vee \cap \h^{-1}
{I_{\scriptscriptstyle F_\h}}^{\hskip-3pt 2} = J_\h + J^\vee
J_\h \; $.  Moreover, multiplication by  $ \h^{-1} $  yields
a  $ \kh $--module  isomorphism  $ \, \mu \, \colon \, J \,
{\buildrel \cong \over
{\lhook\joinrel\relbar\joinrel\relbar\joinrel\twoheadrightarrow}}
\, J^\vee $.  Let  $ \; \rho \, \colon \, J_{\scriptscriptstyle G}
\relbar\joinrel\twoheadrightarrow J_{\scriptscriptstyle G} \big/
{J_{\scriptscriptstyle G}}^{\!2} = \gerg^* \, $  be the natural
projection map, and  $ \; \nu \, \colon \, \gerg^*
\lhook\joinrel\longrightarrow J_{\scriptscriptstyle G} \, $
a section of  $ \rho \, $.  The specialization map  $ \; \pi \,
\colon \, F_\h \relbar\joinrel\twoheadrightarrow
F_0 \; $  restricts to  $ \; \pi' \colon \, J
\relbar\joinrel\twoheadrightarrow J \big/ (J \cap \h \, F_\h)
= J_\h \big/ \h \, J_\h = J_{\scriptscriptstyle G} \, $:  \,
we fix a section  $ \, \gamma \, \colon \, J_{\scriptscriptstyle G}
\lhook\joinrel\relbar\joinrel\rightarrow J_\h \, $  of  $ \pi' $.
Then the composition map  $ \, \sigma := \eta \circ \mu \circ
\gamma \circ \nu \, \colon \, \gerg^* \longrightarrow \gert \, $
is a well-defined Lie bialgebra isomorphism, independent of the
choice of  $ \nu $  and  $ \gamma \, $.  In fact, one has (see [Ga1])
$ \, \fhg \cong \big(\Bbbk[[j_1,\dots,j_n]]\big)[[\h]] \, $  and  $ \,
\uhg \cong \big(\Bbbk\big[j_1^\vee,\dots,j_n^\vee\big]\big)[[\h]] \, $
as topological  $ \kh $--modules.
                                          \par
   For our purposes we need a special choice of the  $ \Bbbk $--basis
$ \, \{y_1, \dots, y_n\} \, $  of  $ \, \gerg^* = J_{\scriptscriptstyle G}
\Big/ {J_{\scriptscriptstyle G}}^{\!2} \, $.  Namely, letting  $ \, k :=
\dim(K) \, $,  \, we fix a system of parameters  $ \{ j_1, \dots, j_k,
j_{k+1}, \dots, j_n \} $  for  $ F[[G]] $  like in the end of \S
1.6: then in particular  $ \, \big(\, \{ j_{k+1}, \dots, j_n \} \!
\mod {J_{\scriptscriptstyle G}}^{\,2} \,\big) \! \mod \gerk^* \, $  is
a  $ \Bbbk $--basis  of  $ \gerg^* \big/ \gerk^* = \gerk^\perp \, $,
\, the cotangent space of  $ G\big/K $  at the point  $ eK \, $.
                                          \par
   By construction  $ \, \big( \calI + J_{\scriptscriptstyle G}^{\,2} \,
\big) \, \cap \, \hbox{\it Span}\,\big( \{j_1,\dots,j_k\} \big) = \{0\}
\, $  and  $ \, \rho(\calI\,) = \big( \calI + {J_{\scriptscriptstyle
G}}^{\!2} \,\big) \! \mod {{J_{\scriptscriptstyle G}}^{\!2}} = \,
\hbox{\it Span}\,\big( \{y_{k+1},\dots,y_n\} \big) = \gerk^\perp \, $.
Thus we choose  {\sl this\/}  set  $ \, \big\{ y_1, \dots, y_k, y_{k+1},
\dots, y_n \big\} \, $  as the basis of  $ \, J_{\scriptscriptstyle G}
\Big/ {J_{\scriptscriptstyle G}}^{\!2} = \gerg^* \, $  to start with.
Then  $ \calI_\h $  identifies with the left ideal of  $ \, \fhg =
\big(\Bbbk[[j_1,\dots,j_n]]\big)[[\h]] \, $  generated by  $ \, \{
j_{k+1},\dots,j_n \} \, $,  \, which is the set of all formal power
series in  $ \, \{ j_1, \dots, j_n, \h \} \, $  such that in each monomial
with non-zero coefficient at least one out of  $ j_{k+1} $,  $ \dots $,
$ j_n $  does occur with non-zero exponent.  Similarly,  $ {\calI_\h}^{\!
\curlyvee} $  identifies with the left ideal of  $ \, \uhg = \big( \Bbbk
\big[j_1^\vee,\dots,j_n^\vee\big] \big)[[\h]] \, $  generated by  $ \,
\big\{ j_{k+1}^\vee,\dots,j_n^\vee \big\} \, $,  \, which is the set of
all formal power series in  $ \h $  with coefficients in  $ \, \Bbbk
\big[j_1^\vee,\dots,j_n^\vee\big] \, $  such that in each monomial in
the  $ j_r^\vee $'s  with non-zero coefficient at least one out of
$ j_{k+1}^\vee $,  $ \dots $,  $ j_n^\vee $  occurs with non-zero
exponent.  But then it's clear   --- thanks to  $ (2.3)' $  ---   that
$ \, {\calI_\h}^{\!\curlyvee} \cap \h \, {F_\h[G]}^\vee \subseteq \,
\h \, {\calI_\h}^{\!\curlyvee} \, $.  The converse inclusion  $ \,
{\calI_\h}^{\!\curlyvee} \cap \h \, {F_\h[G]}^\vee \supseteq \,
\h \, {\calI_\h}^{\!\curlyvee} \, $  is obvious.  Similarly one
proves  {\it (b)}.   \qed
\enddemo

\vskip3pt

   Altogether, Lemmas 4.1--4.3 yield the main result of this section,
namely

\vskip7pt

\proclaim{Proposition 4.4}  $ {\calI_\h}^{\!\curlyvee} \! $,  $ \,
{\gerC_\h}^{\!\! \Lsh} $,  $ \, {\calC_\h}^{\!\!\triangledown} $
and  $ \, {\gerI_\h}^{\! !} $  are quantizations of a pair  {\sl
(coisotropic subgroup, Poisson quotient)},  in the sense of \S 2.7,
for the dual Poisson group  $ G^* $.   \qed
\endproclaim

\vskip3pt

   Next result instead shows that the construction by Drinfeld-like
functors is involutive:

\vskip7pt

\proclaim{Proposition 4.5}  The following identities hold:
 \vskip-11pt
  $$  {\big( \calI_\h^{\,\curlyvee} \big)}^! \, = \; \calI_\h \quad ,
\qquad  {\big( \calC_\h^{\,\triangledown} \big)}^\Lsh \, = \; \calC_\h
\quad ,  \qquad  {\big( {\gerI_\h}^{\! ! \,} \big)}^{\!\curlyvee} \, =
\; \gerI_\h \quad ,  \qquad  {\big( {\gerC_\h}^{\!\!\Lsh\,} \big)}^{\!
\triangledown} \, = \; \gerC_\h \quad .  $$
 \vskip-3pt
\endproclaim

\demo{Proof}  From the very definitions we get
  $$  \displaylines{
   \quad  \delta_n\big(\calI_\h\big) \; \subseteq \; {\textstyle
\sum_{s=1}^n} {J_{F_\h}}^{\hskip-3pt \widetilde\otimes\, (s-1)}
\otimestilde \calI_\h \otimestilde {J_{F_\h}}^{\hskip-3pt
\widetilde\otimes\, (n-s)} \; \subseteq \; {\textstyle \sum_{s=1}^n}
\Big( \h^{s-1} {\big( {F_\h}^{\!\vee} \big)}^{\hskip-1pt \widehat\otimes\,
(s-1)} \Big) \otimeshat   \hfill  \cr
   \hfill   \otimeshat \Big( \h \, \calI_\h^{\,\curlyvee} \Big) \otimeshat
\Big( \h^{n-s} {\big( {F_\h}^{\!\vee} \big)}^{\hskip-1pt \widehat\otimes\,
(n-s)} \Big) \; = \; \h^n \cdot {\textstyle \sum_{s=1}^n} {\big(
{F_\h}^{\!\vee} \big)}^{\hskip-1pt \widehat\otimes\, (s-1)} \otimeshat
\calI_\h^{\,\curlyvee} \otimeshat {\big( {F_\h}^{\!\vee} \big)}^{\hskip-1pt
\widehat\otimes\, (n-s)}  \cr }  $$
for all  $ \, n \in \N_+ \, $,  \, which means exactly that  $ \,
{\big( \calI_\h^{\,\curlyvee} \big)}^! \supseteq \calI_\h \, $.
Similarly, we have also  $ \; \delta_n\big(\calC_\h\big) \, \subseteq \,
{J_{F_\h}}^{\hskip-3pt \widetilde\otimes\, (n-1)} \otimestilde \calC_\h \,
\subseteq \, \Big( \! \h^{n-1} {\big( {F_\h}^{\!\vee} \big)}^{\hskip-1pt
\widehat\otimes\, (n-1)} \Big) \otimeshat \Big( \! \h \, \calC_\h^{\,
\triangledown} \Big) \, = \, \h^n \cdot {\big( {F_\h}^{\!\vee}
\big)}^{\hskip-1pt \widehat\otimes\, (n-1)} \otimeshat \calC_\h^{\,
\triangledown} \; $
for all  $ \, n \in \N_+ \, $,  \, which means exactly that  $ \,
{\big( \calC_\h^{\,\triangledown} \big)}^\Lsh \supseteq \calC_\h \, $.
On the other hand, by definitions  $ \; {\gerI_\h}^{\! ! \,} \cap \,
J_{F_\h} = \, \epsilon \big(\, {\gerI_\h}^{\! ! \,} \cap \, J_{F_\h}
\big) + \, \delta_1 \big(\, {\gerI_\h}^{\! ! \,} \cap \, J_{F_\h} \big)
= \, \delta_1 \big(\, {\gerI_\h}^{\! ! \,} \cap \, J_{F_\h} \big)
\subseteq \h \, \big(\, \gerI_\h \cap \, J_{F_\h} \big) \, $,  \;
which implies  $ \; {\big( {\gerI_\h}^{\! ! \,} \big)}^{\!\curlyvee}
\subseteq \, \gerI_\h \; $.  Similarly,  $ \; {\gerC_\h}^{\!\!\Lsh\,}
\! \cap J_{F_\h} = \, \epsilon \big(\, {\gerC_\h}^{\!\!\Lsh\,} \cap \,
J_{F_\h} \big) + \, \delta_1 \big(\, {\gerC_\h}^{\!\!\Lsh\,} \cap \,
J_{F_\h} \big) = \, \delta_1 \big(\, {\gerC_\h}^{\!\!\Lsh\,} \cap \,
J_{F_\h} \big) \subseteq \, \h \cdot \big(\, \gerC_\h \cap \, J_{F_\h}
\big) \; $  yields  $ \; {\big( {\gerC_\h}^{\!\!\Lsh\,} \big)}^{\!
\triangledown} \subseteq \, \gerC_\h \; $.  Thus all identities in the
claim are half proved.
                                             \par
   To prove the reverse inclusions  $ \; {\big( \calI_\h^{\,\curlyvee}
\big)}^! \subseteq \calI_\h \; $  and  $ \; {\big( \calC_\h^{\,
\triangledown} \big)}^\Lsh \subseteq \, \calC_\h \; $  one can resume
the proof of Proposition 3.2 in [Ga1], which shows that  $ \, {\big(
{F_\h}^{\!\vee}\big)}' \subseteq F_\h \, $:  \; in fact, the same
arguments apply almost untouched with  $ \calC_\h $  instead of
$ F_\h \, $,  and also (with minimal changes) with  $ \calI_\h $
instead of  $ F_\h \, $.  The outcome is  $ \, {\big( \calI_\h^{\,
\curlyvee} \big)}^! \subseteq \calI_\h \, $  and  $ \, {\big(
\calC_\h^{\,\triangledown} \big)}^\Lsh \subseteq \, \calC_\h \, $,
\; whence identities hold.
                                             \par
   To finish with, by Proposition 4.4 we can apply twice Lemma 4.1
and get  $ \, {\big( {\gerC_\h}^{\!\! \Lsh\,} \big)}^{\!\triangledown}
\! = \Big(\! {\big( \calI_\h^{\,\curlyvee} \big)}^{!\,} \Big)^\perp \, $ 
and  $ \; {\big( {\gerI_\h}^{\! !\,} \big)}^{\!\curlyvee} \! = \Big(\!
{\big( \calC_\h^{\,\triangledown} \big)}^{\!\Lsh} \Big)^\perp $.  As 
$ \, {\big( \calI_\h^{\, \curlyvee} \big)}^! \! = \, \calI_\h \, $
and  $ \, {\big( \calC_\h^{\, \triangledown} \big)}^\Lsh \! = \, \calC_\h
\, $,  \, we get  $ \; {\big( {\gerC_\h}^{\!\!\Lsh\,} \big)}^{\!
\triangledown} \! = {\calI_\h}^{\!\perp} \, $  and  $ \, {\big(
{\gerI_\h}^{\! !\,} \big)}^{\!\curlyvee} \! = {\calC_\h}^{\!\perp} \, $;
\, but then (2.4) eventually yields  $ \, {\big( {\gerC_\h}^{\!\!\Lsh\,}
\big)}^{\!\triangledown} \! = \, \gerC_\h \, $  and  $ \, {\big(
{\gerI_\h}^{\! !\,} \big)}^{\!\curlyvee} \! = \, \gerI_\h \, $.   \qed
\enddemo

 \vskip4pt

   {\it  $ \underline{\hbox{\it Remark}} $:} \, like for Lemma 4.2,
in the previous proof we might prove only one of the identities in
the claim, e.g.~that for  $ \calI_\h \, $:  \, all others then follow
via  {\it (i)--(iv)\/}  in Lemma 4.1.

\vskip1,3truecm

\centerline {\bf \S\; 5 \ Specialization at  $ \, \h = 0 \, $ }

\vskip11pt

   We shall now look at semiclassical limits of the images of
Drinfeld-like functors.  The result   --- Proposition 5.2 ---
will be  $ \, \Big( K^\perp, \, G^* \big/ K^\perp \Big) \, $,
\, in the sense that this will be the pair  {\sl (coisotropic
subgroup, Poisson quotient)\/}  mentioned in Proposition 4.4.

\vskip7pt

\proclaim{Lemma 5.1} Let  $ \Cal{S}(G^*) $  be the set of formal
subgroups of the formal Poisson group  $ G^* \, $.
                             \hfill\break
   \indent   (a) \;  $ {\calI_0}^{\!\curlyvee} \! \ideal_\ell \!
\coideal \, {F_0[[G]]}^\vee \! = \, U(\gerg^*) \, $,  \hskip2pt
   \hbox{whence  \hskip3pt  $ {\calI_0}^{\!\curlyvee} \! = U(\gerg^*)
\! \cdot \! \frak{l} $  \hskip2pt  for some Lie subalgebra  \hskip1pt
$ \frak{l} \leq \gerg^* $;}
%
   \indent   (b) \;  $ {\calC_0}^{\!\!\triangledown} \leq_\Hpicc
{F_0[[G]]}^\vee = \, U(\gerg^*) \, $,  \hskip2pt  whence  \hskip2pt
$ \; {\calC_0}^{\!\!\triangledown} = U(\gerh) \; $  \hskip2pt  for
some Lie subalgebra  \hskip2pt  $ \, \gerh \leq \gerg^* \, $;
                              \hfill\break
   \indent   (c) \;  $ {\gerI_0}^{\! !} \ideal_\Hpicc {U_0(\gerg)}' = \,
F[[G^*]] \, $,  \hskip2pt  whence  \hskip2pt  $ \; {\gerI_0}^{\! !}
= \calI(\varGamma) \; $  \hskip2pt  for some  $ \, \varGamma \in
\Cal{S}(G^*) \, $;
                             \hfill\break
   \indent   (d) \;  $ {\gerC_0}^{\!\!\Lsh} \! \leq^1 \!\! \coideal_\ell\,
{U_0(\gerg)}' = \, F[[G^*]] \, $,  \hskip1pt  whence  \hskip3pt
$ \; {\gerC_0}^{\!\!\Lsh} \! = \! {F[[G^*]]}^\varTheta \; $  \hskip2pt
for some  \hskip2pt  $ \, \varTheta \in \Cal{S}(G^*) \, $;
                             \hfill\break
   \indent   (e) \,  Let  $ \, H \in \Cal{S}(G^*) \, $  be the formal
subgroup of  $ \, G^* \, $  with  $ \, \text{\it Lie}\,(H) = \gerh \, $,  and
let  $ \, L \in \Cal{S}(G^*) \, $  be the one with  $ \, \text{\it Lie}\,(L) = \frak{l} \, $.  Then  $ \, \varGamma = H = L = \varTheta \, $.
                             \hfill\break
   \indent   (f) \,  the formal subgroup  $ \, \varGamma = H = L = \varTheta
\, $  in (e) is coisotropic in  $ G^* $.
\endproclaim

\demo{Proof}  Statements  {\it (a)\/  {\rm and}  (d)\/}  follow
trivially from Lemma 4.2; the same also implies part of  {\it (b)\/
{\rm and}  (c)},  in that  $ {\gerI_0}^{\! !} $  is a  {\sl bialgebra
ideal\/}  of  $ {U_0(\gerg)}' $  and  $ {\calC_0}^{\!\!\triangledown} $
is a  {\sl subbialgebra\/}  of  $ {F_0[[G]]}^\vee \, $.  Now,  $ \,
{F_0[[G]]}^\vee = U(\gerg^*) \, $,  \, and a subbialgebra of any
universal enveloping algebra  (such as  $ U(\gerg^*) \, $)  is
automatically a Hopf subalgebra: thus $ {\calC_0}^{\!\!\triangledown} $
is a Hopf subalgebra.  On the other hand, the orthogonality relations
of  Lemma 5.1{\it (ii)\/}  imply that  $ {\gerI_0}^{\! !} $  is a
Hopf ideal too.
                                                    \par
  Claim  {\it (e)\/}  follows directly from Proposition 4.4 and from
Remark 2.8{\it (b)}.
                                                    \par
   Finally  {\it (f)\/}  follows from Proposition 4.4 and
Lemma 2.9.  \qed
\enddemo

\vskip7pt

\proclaim{Proposition 5.2} The coisotropic subgroup  $  \, \varGamma = H
= L = \varTheta \, $  of Proposition 5.1 coincide with  $ \, K^\perp \in
\hbox{\it co}\hskip1pt\Cal{S}(G^*) \, $  (cf.~Definition 1.4).  In other
words,  $ \, \gerl = \gerh \, $  coincides with  $ \, \gerk^\perp \, \big(
\subseteq \gerg^* \big) \, $.
\endproclaim

\demo{Proof} We resume the construction made for the proof of Lemma 4.3,
with same notation.  In particular we fix a special subset  $ \, \{ j_1,
\dots, j_k, j_{k+1}, \dots, j_n \} \, $  of  $ J_{\scriptscriptstyle G} $
enjoying the properties mentioned there, and call  $ \, \{ y_1, \dots,
y_k, y_{k+1}, \dots, y_n \} \, $  its image in  $ \, \gerg^* =
J_{\scriptscriptstyle G} \big/ {J_{\scriptscriptstyle G}}^{\! 2} \, $.
                                          \par
%
%
   The same kind of analysis carried on in the proof of Lemma 4.3
to prove that  $ \, \sigma \, \colon \, \gerg^* \! \cong \gert \, $
shows that
%
%
the unital subalgebra  $ \, {\calC_0}^{\!\!\triangledown}
:= {\calC_\h}^{\!\!\triangledown} \!\! \mod \h \, {F_\h}^{\!\vee} \, $
is generated by  $ \, \eta \big( {\calC_\h}^{\!\!\triangledown} \cap
J^\vee \big) = (\mu \circ \eta) \big( \calC_\h \cap J \,\big) =
(\sigma \circ \rho \circ \pi) \big( \calC_\h \cap J \,\big) =
\sigma \big( \rho( \calC \cap J_{\scriptscriptstyle G}) \big) =
\sigma \big( \rho( \langle j_{k+1}, \dots, j_n \rangle ) \big) =
\sigma \big( \gerk^\perp \big) \, $,  \, where  $ \, \langle\, j_{k+1},
\dots, j_n \rangle \, $  is the ideal of  $ \calC $  generated by
$ \, \{j_{k+1}, \dots, j_n\} \, $.  Therefore  $ \, {\calC_0}^{\!
\!\triangledown} = U(\gerh) \, $  is generated by  $ \gerk^\perp $,
whose elements are primitive, so belong to  $ \gerh \, $:  \,
then  $ \, \gerh = \gerk^\perp \, $,   \, q.e.d.   \qed
\enddemo

\vskip7pt

\proclaim{Corollary 5.3}  $ \, \calI(K)_\h^{\,\curlyvee} \, $,  $ \,
\calC(K)_\h^{\,\triangledown} \, $,  $ \, \gerI(K)_\h^{\; !} \, $
and  $ \, \gerC(K)_\h^{\,\Lsh} \, $  all provide quantizations,
{w.r.t.} $ \big( {U_\h}^{\!\prime}, {F_\h}^{\!\vee} \big) \, $,
of the formal coisotropic subgroup  $ K^\perp $  and the
formal Poisson quotient  $ \, G^* \! \big/ K^\perp $.
\endproclaim

\demo{Proof}  The claim follows from Proposition 4.4, Lemma 5.1
and Proposition 5.2.   \qed
\enddemo

\vskip5pt

  Patching together all previous results, we can finally prove
Theorem 3.3:

\vskip7pt

  $ \underline{\hbox{\it Proof of Theorem 3.3.}} \, $  Corollary 5.3
proves that the functors in  {\it (a)\/}  are well-defined on objects,
and it is trivially clear that they are inclusion-preserving, so they
do are functors.  Proposition 4.5 proves the rest of claim  {\it (a)},
in particular that these functors are in fact equivalences.  In addition,
Corollary 5.3 also proves claim  {\it (b)}.   \qed

 \vskip1,3truecm

\centerline {\bf \S\; 6 \ Example: the Stokes matrices as Poisson
homogeneous  $ {{SL}_n}^{\!\!*} $--space }

\vskip11pt

  {\bf 6.1 The Poisson homogeneous  $ {{SL}_n}^{\!\!*} $--space  of Stokes
matrices.} \, Let  $ \, G = SL_n(\Bbbk) \, $  endowed with the standard
Poisson-Lie structure.  We denote by  $ \frak{d} $  the Cartan
subalgebra of diagonal matrices in  $ \gersln(\Bbbk) \, $.  With
$ \frak{b}_+ $ (resp.~$ \frak{b}_- $)  we denote the Borel
        \hbox{subalgebra of upper}
 \eject
\noindent
 (resp.{} lower) triangular matrices in  $ \gersln $;  then
$ B_+ $  and  $ B_- $  will be the corresponding Borel subgroups in
$ SL_n \, $.  It is well known that at the infinitesimal level the
dual Lie bialgebra can be identified with  $ \; \gerg^* \, = \, \big\{
(X,Y) \in \frak{b}_+ \oplus \frak{b_-} \,\big|\, X\big|_{\frak{d}} =
-Y\big|_{\frak{d}} \big\} \, $,  \; so that the simply connected dual
Poisson group is  $ \, G^* = B_+ \star B_- \, $,  \, the pairs of upper
and lower triangular matrices such that the restrictions on the diagonal
are mutually inverse.
                                          \par
   By construction, the algebra  $ \, F[G^*] = F[B_+ \star B_-] \, $
is generated by matrix coefficients  $ x_{i,j} $  ($ 1 \leq i \leq j
\leq n $)  for the over-diagonal part of  $ B_+ \, $,  $ y_{i,j} $
($ 1 \geq i \geq j \geq 1 $)  for the under-diagonal part of
$ B_- \, $,  and  $ z_i $  ($ 1 \leq i \leq n $)  for the
diagonal part of  $ B_+ \, $.
                                        \par
   Let  $ \, H = SO_n(\Bbbk) \lhook\joinrel\relbar\joinrel\longrightarrow
SL_n(\Bbbk) \, $  be the standard embedding.  The corresponding Lie algebra
is  $ \, \gerh = \gerso_n(\Bbbk) \, $.  Its orthogonal in  $ \gerg^* $,
for the pairing given by the Killing form, is  $ \; \gerh^\perp \, =
\, \big\{ \big( b, -b^{\,t} \big) \in \frak{b}_+ \oplus \frak{b}_-
\,:\, b\big|_{\frak{d}} = 0 \big\} \; $  and can be integrated to
$ \; H^\perp \, = \, \big\{ (B,C) \in B_+ \star B_- \,\big|\, B \,
C^{\,t} = \text{\it Id} \,\big\} \, $,  \; which is a coisotropic
subgroup of  $ G^* \, $.  We are then in the situation described in
\S 1.  The spaces  $ \, SL_n \big/ SO_n \, $  and  $ \, {{SL}_n}^{\!\!*}
\big/ H^\perp \, $  are a complementary dual pair of Poisson homogeneous
spaces: the former can be identified with the space of symmetric matrices
and the latter with the space  $ U^+_n $  of Stokes matrices, i.e.~upper
triangular unipotent  $ (n \times n) $--matrices.  By construction the
function algebra  $ \, F\big[U^+_n\big] \! = F\big[G^*\big/H^\perp\big]
\! = {F[G^*]}^{H^\perp} $  is generated by elements  $ \, x_{i,j} \, $,
\, for all  $ \, 1 \leq i < j \leq n \, $,  \, which may be realized
as the matrix coefficient functions on Stokes matrices.
                                          \par
   The Poisson structure on  $ U^+_n $  was first found by Dubrovin
in the  $ \, n = 3 \, $  case (see [Du]) and then by Ugaglia (cf.~[Ug])
for generic  $ \, n \geq 3 \, $  in a completely different setting: it
naturally arises in the study of moduli spaces of semisimple Frobenius
manifolds.  Later, in [Bo,Xu], it was shown how  $ U^+_n $  with such
structure is a Poisson homogeneous space of the Poisson-Lie group
$ \, B_+ \star B_- \, $,  \, dual to the standard  $ SL_n $,  as just
explained.  More explicitly, from [Xu] one can argue the following

\vskip7pt

\proclaim{Proposition 6.2} \, Let  $ \; \Psi : B_+ \star B_- \!
\loongrightarrow B_+ \star B_- \, $,  $ \, \Psi(B,C) := \big( C^{\,t},
B^{\,t} \big) \, $  and let  $ \, H^\perp = \big\{\, g \in B_+ \star
B_- \;\big|\; \psi(g) = g^{-1} \,\big\} \, $.  Then  $ H^\perp $  is
a coisotropic subgroup of  $ \, B_+ \star B_- $  and  $ \, U^+_n \cong
(B_+ \star B_-) \big/ H^\perp \, $  with its quotient Poisson structure.
\qed
\endproclaim

\vskip3pt

  {\bf 6.3 Towards quantization of Stokes matrices.} \, In the present
section we look for quantizations of  $ U^+_n \, $:  \, the first step
is to switch to the associated formal homogeneous space.  Actually,
the function algebra  $ \, F_\h\big[\big[G^*\big/H^\perp\big]\big] =
F_\h\big[\big[U^+_n\big]\big] \, $  is nothing but the algebra of formal
power series in the matrix coefficient functions, say  $ \, \chi_{i,j}
\, $  ($ \, 1 \leq i < j \leq n \, $),  on  $ \, U^+_n \, $.
                                                  \par
   Now we look for a quantization  $ \, F_\h\big[\big[U^+_n\big]\big] \, $
of  $ F \big[\big[U^+_n\big]\big] $  with the above Poisson structure:
we shall find it applying Theorem 3.3.  As our purpose is to obtain a
quantum algebra of functions on the homogeneous space, an object of type
{\it (b)\/}  in the list (2.3), we start with an object of type  {\it
(d)\/}  in the same list.  This means that as a starting point we need
a subalgebra and left coideal inside  $ U_\h(\gersln) $  quantizing the
standard embedding of  $ \gerso_n \, $.  This has been already obtained
in [No, \S 2.3] (see also the works of Klimyk et al., e.g~[GIK] and
references therein): we recall hereafter its definition in the formal
setup.  We begin fixing notation for  $ U_\h(\gergln) $,  a quantum
analogue of  $ \, U(\gergln) \, $,  \, and its Hopf subalgebra
$ U_\h(\gersln) \, $:

\vskip7pt

\proclaim{Definition 6.4} \, We call  $ U_\h(\gergln) $  the topological,
$ \h $--adically  complete, associative unital  $ \kh $--algebra  with
generators  $ f_i \, $,  $ \ell_j \, $,  $ e_i \, $  ($ \, i = 1, \dots,
n-1 $;  $ j = 1, \dots, n $)  and relations
 \vskip-11pt
  $$  \displaylines{
   \hfill   \ell_j \, f_i - f_i \, \ell_j = \, (\delta_{i+1,j} -
\delta_{i,j}) \, f_i \; ,  \quad  \ell_j \, \ell_k = \ell_k \, \ell_j \; ,
\quad  \ell_j \, e_i - e_i \, \ell_j = \, (\delta_{i,j} - \delta_{i+1,j})
\, e_i  \quad   \hfill  \forall \;\; i \, , j \, , k  \cr
 }  $$
 \eject
  $$  \displaylines{
   \hskip1pt   e_k \, f_l \, - \, f_l \, e_k \; = \; \delta_{k,l}
\, {{\; t_k^{+1} - t_k^{-1} \;} \over {\, q - q^{-1} \,}}
\quad  \forall \;\; k \, , l \, ,  \qquad
  e_i \, e_j \; = \; e_j \, e_i \; ,  \,\;\;
f_i \, f_j \; = \; f_j \, f_i   \quad
\forall \;\; \vert\, i \!-\! j \,\vert > 1  \cr
   \hfill   e_i^2 \, e_j - \left( q + q^{-1} \right) e_i \, e_j \, e_i
+ e_j \, e_i^2 \, = \, 0 \; ,  \;  f_i^2 \, f_j - \left( q + q^{-1}
\right) f_i \, f_j \, f_i + f_j \, f_i^2 \, = \, 0   \hfill  \;
\forall \; \vert\, i \!-\! j \,\vert = 1  \cr }  $$
where hereafter we use notation  $ \, q:= \exp(\h) \, $,  $ \, q^X :=
\exp\big(\h \, X\big) \, $  and  $ \, t_i := q^{\ell_i - \ell_{i+1}}
\, $  ($ \, \forall \, i \, $).  It has a structure of topological
Hopf  $ \kh $--algebra  uniquely given by
 \vskip-11pt
  $$  \displaylines{
   \hfill   \Delta(f_i) \, = \, f_i \otimes t_i^{-1} + 1 \otimes f_i \; ,
\; \qquad  S(f_i) \, = \, - f_i \, t_i \; , \; \qquad  \epsilon(f_i) = 0
\hfill \forall \;\; i \, \phantom{.}  \cr
   \hfill   \qquad  \Delta(\ell_j) \, = \, \ell_j \otimes 1 + 1 \otimes
\ell_j \; ,  \qquad \qquad \qquad  S(\ell_j) \, = \, - \ell_j \; ,
\quad \qquad  \epsilon(\ell_j) \, = \, 0   \hfill  \forall \;\; j
\, \phantom{.}  \cr
   \hfill   \;\; \Delta(e_i) \, = \, e_i \otimes 1 + t_i \otimes e_i
\; ,  \;\;\; \qquad  S(e_i) \, = \, - t_i^{-1} \, e_i \; ,  \; \qquad
\epsilon(e_i) \, = \, 0   \hfill  \forall \;\; i \, \phantom{.}  \cr }  $$
\endproclaim
%
%
%

\vskip5pt

  {\bf 6.5  Quantum root vectors and  $ L $--operators  in
$ U_\h(\gergln) \, $.} \, We recall the notion of  {\sl
$ L $--operators},  first introduced in [FRT]: these are
elements  $ \, L^\pm_{i,j} \in U_\h(\gergln) \, $  (with
$ \, i, j = 1, \dots, n $),  which are defined as follows.
Set  $ \, [x,y]_a := x \, y - a \, y \, x \, $  (for all
$ x $, $ y $,  $ a \, $),  and define
  $$  E_{i,i+1} := e_i \; ,  \quad  E_{i,j} := \big[ E_{i,k},
E_{k,j} \big]_q \; ,  \; \quad  F_{i+1,i} := f_i \; ,  \quad
F_{j,i} := \big[ F_{j,k}, F_{k,i} \big]_{q^{-1}} \;   \eqno
\forall \;  i \! < \! k \! < \! j  $$
(where  $ \, q := \exp(\h) \, $  again).  These are  {\sl quantum root
vectors\/}  in  $ U_\h(\gergln) \, $,  in that the coset of  $ \, E_{i,j}
\, $  (resp.~$ \, F_{j,i} \, $)  modulo  $ \, \h \, U_\h(\gergln) \, $
in  $ \, U_\h(\gergln) \Big/ \h \, U_\h(\gergln) \cong U(\gergln) \, $
is the elementary matrix  $ \text{e}_{i,j} $  (resp.~$ \text{e}_{j,i} $)
for all  $ \, i < j \, $.
                                        \par
   The  $ L $--operators  are obtained by twisting and rescaling the
above quantum root vectors,
  $$  \displaylines{
   {} \hfill   L^+_{i,i} := \, q^{+\ell_i} =: g_i^{+1} \; ,  \qquad
L^+_{i,j} := + \big( q - q^{-1} \big) \, g_i^{+1} F_{j,i} \;\, ,
\qquad  L^+_{j,i} := \, 0  \hfill   (i<j)  \cr
   {} \hfill   L^-_{i,i} := \, q^{-\ell_i} =: g_i^{-1} \; ,  \qquad
L^-_{i,j} := - \big( q - q^{-1} \big) \, E_{j,i} \, g_j^{-1} \; ,
\qquad  L^-_{j,i} := \, 0  \hfill   (i>j)  \cr }  $$
and satisfy the remarkable formulas  $ \; \Delta \big( L^\pm_{i,j}
\big) \, = \, {\textstyle \sum_{k = i \wedge j}^{i \vee j}} \,
L^\pm_{(i \wedge j), k} \otimes L^\pm_{k, (i \vee j)} \; $,  $ \;
\epsilon\big(L^\pm_{i,j}\big) \, = \, \delta_{i,j} \; $.
                                               \par
   When suitably normalized, the  $ L $--operators  are again
$ q $--analogues  of the elementary matrices of  $ \gergln \, $:
namely, the coset of  $ \, \big( q - q^{-1} \big)^{-1} L^+_{i,j}
\, $  (resp.~$ \, \big( q - q^{-1} \big)^{-1} L^-_{j,i} \, $)  modulo
$ \, \h \, U_\h(\gergln) \, $  in the semiclassical limit  $ \,
U_\h(\gergln) \Big/ \h \, U_\h(\gergln) \cong U(\gergln) \, $  is
$ \text{e}_{j,i} $  (resp.~$ \text{e}_{i,j} $) for all  $ \, i < j
\, $.  Moreover, the elements  $ \, \widehat{L}^\pm_{i,j} := \big(
q - q^{-1} \big)^{\delta_{i,j} - 1} L^\pm_{i,j} \, $  for  $ \, i
\! \not= \! j \, $  together with the  $ \ell_k $'s  form a set of
generators for  $ U_\h(\gergln) \, $.  Set also  $ \, \Lambda^\pm
:= \big( \Lambda^\pm_{i,j} \big)_{i,j=1}^n \, $  for any  $ \,
\Lambda \in \big\{ L, \widehat{L} \,\big\} \, $.

\vskip7pt

   {\bf 6.6 Quantization of  $ U(\gersln) \, $.} \, For all  $ \, i = 1,
\dots, n-1 \, $,  \, let  $ \, h_i := \ell_i - \ell_{i+1} \, $.  Given
$ U_\h(\gergln) $  as above, we define  $ U_\h(\gersln) $  as the closed
topological subalgebra of  $ U_\h(\gergln) $  generated by  $ \, \big\{
f_i \, , h_i \, , e_i \,\big\}_{i = 1, \dots n-1} \, $.  From the
presentation of  $ U_\h(\gergln) $  in Definition 6.4 one argues
a presentation of  $ U_\h(\gersln) $  as well: in particular, this
shows that  $ U_\h(\gersln) $  is a  {\sl Hopf subalgebra\/}  of
$ U_\h(\gergln) \, $;  \, moreover, by construction we have a
quantum analogue of the classical embedding  $ \, \gersln
\! \lhook\joinrel\longrightarrow \gergln \, $.  Note also
that, for any  $ i $,  $ j \, $,  we have  $ \, L^\pm_{i,j}
\not\in U_\h(\gersln) \, $.  It is also immediate to check that
our Hopf algebra  $ U_\h(\gersln) $  coincides with Drinfeld's one.

\vskip7pt

   {\bf 6.7 Quantization of  $ U(\gerso_n) \, $.} \, Following an
idea of Noumi, Klimyk et al., we define  $ U_\h(\gerso_n) $  as a
subalgebra of  $ U_\h(\gersln) $.  We call  $ U_\h(\gerso_n) $
the closed topological  $ \kh $--subalgebra of  $ U_\h(\gergln) $
generated by the matrix entries of  $ \; K \, := \, \big( \widehat{L}^-
\big)^t \, J \, L^+ \, = \, \big( L^- \big)^t \, J \, \widehat{L}^+ $,
\; where  $ J $  is the  $ (n \times n) $  diagonal matrix  $ \,
\text{\sl diag}\,\big( q^{n-1}, \dots, q, 1 \big) \, $.  Explicit
computations give
 \vskip-13pt
  $$  K_{i,j} \, = \, {\textstyle \sum\nolimits_{k=i}^j} \, q^{n-k}
\big( q - q^{-1} \big)^{-1} L^-_{k,i} \, L^+_{k,j} \, = \, {\textstyle
\sum\nolimits_{k=i}^j} \, q^{n-k} \widehat{L}^-_{k,i} \, L^+_{k,j} \,
= \, {\textstyle \sum\nolimits_{k=i}^j} \, q^{n-k} L^-_{k,i} \,
\widehat{L}^+_{k,j}  $$
 \vskip-1pt
\noindent
 for the matrix entries of  $ K $,  which is upper triangular with
$ J $  onto the diagonal.  Note that we have  $ \, \big( q - q^{-1}
\big)^{\delta_{i,j} - 1} L^-_{k,i} \, L^+_{k,j} \in U_\h(\gersln) \, $
for all  $ i, k, j $,  hence  $ \, U_\h(\gerso_n) \subseteq U_\h(\gersln)
\, $  as well.  This yields quantum analogues of the classical
embeddings  $ \, \gerso_n \! \lhook\joinrel\longrightarrow \gersln
\! \lhook\joinrel\longrightarrow \gergln \, $.  Moreover, w.r.t.~the
Lie bialgebra structure on  $ \gerg $  inherited by its quantization
$ U_\h(\gergln) $  one has that  $ \gerso_n $  is also a Lie coideal
of  $ \gergln \, $,  \, hence correspondingly  $ SO_n $  is a coisotropic
subgroup of  $ GL_n \, $.  Note that we have fixed Noumi's parameters
$ a_j $  to be  $ \, a_j = q^{n-j} \, $  (for all  $ j \, $).  With
respect to the coproduct,  $ U_\h(\gerso_n) $  {\sl is a right coideal\/}
both of  $ U_\h(\gersln) $  and of  $ U_\h(\gergln) \, $.  Thus  $ \,
\gerC_\h := U_\h(\gerso_n) \, $  and  $ \, U_\h(\gerg) := U_\h(\gersln)
\, $  do realize the situation of  (2.3--({\it d\/}))   --- the
specialization result  $ \, U_\h(\gerso_n)\big|_{\h=0} \! \cong
U(\gerso_n) \, $  being explained in [No] ---   but for having a
{\sl right\/}  instead than  {\sl left\/}  coideal.  However,  {\sl
by left-right symmetry our analysis remains unchanged}.  So  $ \,
\gerC_\h := U_\h(\gerso_n) \, $  is a quantum subgroup for the
quantum group  $ U_\h(\gergln) \, $.
                                                \par
   We now apply the functor  $ \, {(\ )}^\Lsh \, \colon \, \gerC_\h
\big( \hbox{\it co}\hskip1pt\Cal{S}(SL_n) \big) \, {\buildrel \cong
\over \llongrightarrow} \, \calC_\h \big( \hbox{\it co} \hskip1pt
\Cal{S}(B_+ \star B_-) \big) \, $  of Theorem 3.3 to get a quantization 
$ \, F_\h\big[\big[U^+_n\big]\big] := {U_\h(\gerso_n)}^{\!\Lsh} \, $  of 
$ \, F\big[\big[SO_n^\perp\big]\big] = F\big[\big[U^+_n\big]\big]
\, $.  We explain in detail the case of  $ \, n = 3 \, $,  \, and then
basing on that we will give a sketch of the general situation.  Note that
the over-diagonal entries of the matrix  $ K $  will provide   --- passing
from  $ U_\h(\gerso_n) $  to  $ \, F_\h \big[\big[ U^+_n \big]\big] :=
{U_\h(\gerso_n)}^{\!\Lsh} \, $  and eventually to the semiclassical limit
of the latter ---   algebra generators of  $ F\big[\big[U^+_n\big]\big] $, 
namely the matrix coefficients of Stokes matrices.

\vskip4pt

   {\it  $ \underline{\hbox{\it Warning}} $:} \, Noumi's definition
of  $ U_\h(\gerso_n) $  is in [No, \S 2.4] (mutatis mutandis).  It is
explained there that one can take as algebra generators of  $ U_\h
(\gerso_n) $  the entries of either one of four different matrices,
given in formula (2.18) in  [{\it loc.~cit.}].  Among these, we choose
$ \; K_0 := \, \big(L^-\big)^t \, Q \, J^{-1} \, L^+ \, $,  \, where
$ J $  is given above and  $ Q $  is the  $ (n \times n) $  diagonal
matrix  $ \, \text{\sl diag}\,\big( q^{n-1}, \dots, q, 1 \big) = J^2
\, $,  \, so that  $ \, Q \, J^{-1} = J \, $.  We also need to rescale
such generators, and eventually take  $ \, K := \big( q - q^{-1}
\big)^{-1} \, K_0 \, $  as above for the purpose of specialization.

\vskip7pt

   {\bf 6.8 The algebras  $ {U_\h(\gergln)}' $  and  $ {U_\h(\gersln)}'
\, $.} \, As  $ \, F_\h\big[\big[U^+_n\big]\big] := {U_\h(\gerso_n)}^{\!
\Lsh} \, $  is a subalgebra of  $ {U_\h(\gergln)}' $  and  $ {U_\h
(\gersln)}' \, $,  we do need a clear description of these objects.
                                                \par
   By definition, the topological Hopf algebra  $ U_\h(\gergln) $  is
{\sl  $ Q $--graded},  $ Q $  being the root lattice of  $ \gergln \, $,
\, with  $ \, \partial(f_i) = -\alpha_i \, $,  $ \, \partial(h_i) = 0
\, $,  $ \, \partial(e_i) := +\alpha_i \, $  where  $ \, \alpha_i \, $
is the  $ i $--th  simple root of  $ \gergln \, $,  \, for all  $ i \, $.
Also,  $ \, \partial \big( F_{j,i} \big) = \partial \big( \Lambda^+_{i,j}
\big) = - \sum_{i \leq k \leq j} \alpha_k =: - \alpha_{i,j} \, $  and
$ \, \partial \big( E_{i,j} \big) = \partial\big(\Lambda^-_{j,i}\big)
= + \sum_{i \leq k \leq j} \alpha_k =: +\alpha_{i,j} \, $,  \, for
all  $ \, i < j \, $  and  $ \, \Lambda \in  \big\{ L, \widehat{L}
\,\big\} \, $.  It follows that  $ {U_\h(\gergln)}^{\otimes d} $  is
$ Q^{\oplus d} $--\,graded  as a topological algebra, and the like
for  $ {U_\h(\gersln)}^{\otimes d} \, $  (for all $ \, d \in \N \, $).
                                                \par
   The formulas for the coproduct of  $ L $--operators  in \S 6.5
can be iterated, yielding for  $ \widehat{L}^\pm $
  $$  \Delta^d\big(\,\widehat{L}^+_{i,j}\big) \; = \; {\textstyle
\sum_{I^+_d}} \big( q - q^{-1} \big)^{(d - 1 - \delta_{i, k_1} -
\delta_{k_1,k_2} - \cdots - \delta_{k_{d-1}, j})} \, \widehat{L}^+_{i,
k_1} \! \otimes \widehat{L}^+_{k_1,k_2} \! \otimes \cdots \otimes
\widehat{L}^+_{k_{d-1}, j}   \hskip11pt  (\, d \geq 2)  $$
where  $ \; I^+_d := \big\{\, k_1, \dots, k_{dl-1} \,\big|\, i \leq
k_1 \leq k_2 \leq \cdots \leq k_{d-1} \leq j \,\big\} \; $  for
$ \, i < j \, $,  \, and similarly
  $$  \Delta^d\big(\,\widehat{L}^-_{i,j}\big) \; = \; {\textstyle
\sum_{I^-_d}} \big( q - q^{-1} \big)^{(d - 1 - \delta_{i, k_1}
- \delta_{k_1,k_2} - \cdots - \delta_{k_{d-1}, j})} \,
\widehat{L}^-_{i, k_1} \! \otimes \widehat{L}^-_{k_1,k_2}
\! \otimes \cdots \otimes \widehat{L}^-_{k_{d-1}, j}  $$
where  $ \; I^-_d := \big\{\, k_1, \dots, k_{d-1} \,\big|\, i
\geq k_1 \geq k_2 \geq \cdots \geq k_{d-1} \geq j \,\big\} \; $
for  $ \, i > j \, $.  In particular,
  $$  \Delta^d\big(\,\widehat{L}^\varepsilon_{i,j}\big) \; =
\; {\textstyle \sum_{r + s = \ell - 1}} \big( g_i^{\varepsilon 1}
\big)^{\otimes r} \! \otimes \widehat{L}^\varepsilon_{i,j} \otimes
\big( g_j^{\varepsilon 1} \big)^{\otimes s} + \, R   \eqno
(\hbox{hereafter\ } \; \varepsilon \in \{+,-\} )  $$
where  $ R $  is a topological sum of homogeneous terms in  $ {U_\h
(\gergln)}^{\otimes \ell} $  whose degree in  $ Q^{\oplus d} $  is of
type  $ \, (\partial_1, \dots, \partial_d) \, $,  \, each  $ \partial_k $
being a positive or negative root (according to  $ \, \varepsilon = - \, $
or  $ \, \varepsilon = + \, $)  of height less than that of  $ \alpha_{i,j}
\, $.  Finally, for all  $ \, i = 1, \dots, n \, $  we have
  $$  \Delta^d\big(h_i\big) \, = \, {\textstyle \sum_{r+s=d-1}}
1^{\otimes r} \otimes h_i \otimes 1^{\otimes s}   \eqno \forall
\;\; d \in \N_+ \; .  $$
   \indent      Now let  $ \Phi_+ $  (resp.~$ \Phi_- $)  be the
set of positive (resp.~negative) roots of  $ \gergln \, $,  \, and
fix any total ordering  $ \preceq $  on  $ \Phi_+ \, $.  Set also
$ \, L^\pm_\alpha := L^\pm_{i,j} \, $  for each root  $ \, \alpha
= \mp \, \alpha_{i,j} \, $.  The well-known quantum PBW theorem
(adapted to the present case) ensures that
  $$  \Bbb{S} \, := \, \Big\{\, {\textstyle \prod_{\alpha
\in \Phi_-}} \!\! \big(\, \widehat{L}^+_\alpha \big)^{\lambda_\alpha^+}
\, {\textstyle \prod_{i=1}^n} \ell^{\,\eta_i} \, {\textstyle
\prod_{\alpha \in \Phi_+}} \!\! \big(\, \widehat{L}^-_\alpha
\big)^{\lambda_\alpha^-} \;\Big|\; \lambda_\alpha^+ \, , \eta_i \, ,
\lambda_\alpha^+ \in \N \; \forall \, \alpha, i \,\Big\}  $$
is a topological  $ \kh $--basis  of  $ U_\h(\gergln) \, $;  \, hereafter
the products over positive or negative roots are made w.r.t.~the fixed
total ordering.
                                         \par
   Given  $ \, \Cal{M} \in \Bbb{S} \, $  we set  $ \, \big|
\Cal{M} \big| := \sum_{\alpha \in \Phi_-} \lambda_\alpha^+ + \sum_{i=1}^n
\eta_i + \sum_{\alpha \in \Phi_+} \lambda_\alpha^- \, $,  \, the sum of
all exponents occurring in  $ \Cal{M} \, $.  Since  $ \Delta^d $  is
a  {\sl graded algebra morphism},  the previous formulas imply that for
each PBW-like monomial  $ \Cal{M} $  in  $ \Bbb{S} $  we have, for all
$ \, d \geq \big|\Cal{M}\big| \, $,
  $$  \eqalign{
   \Delta^d\big(\Cal{M}\big) \,  &  = \,
\widehat{L}^+_{-\alpha_1} \zeta^{\,(1)}_{-\alpha_1} \otimes \cdots \otimes
\widehat{L}^+_{-\alpha_1} \zeta^{(\lambda^+_{-\alpha_1})}_{-\alpha_1}
\otimes \cdots \otimes \widehat{L}^+_{-\alpha_N} \zeta^{\,(1)}_{-\alpha_N}
\otimes \cdots \otimes \widehat{L}^+_{-\alpha_N}
\zeta^{(\lambda^+_{-\alpha_N})}_{-\alpha_N} \otimes   \hfill  \cr
   {} \otimes h_1 \,  &  \theta^{(1)}_{h_1} \otimes \cdots
\otimes h_1 \theta^{(\eta_1)}_{h_1} \otimes \cdots \otimes
h_{n-1} \theta^{(1)}_{h_{n-1}} \otimes \cdots \otimes h_{n-1}
\theta^{(\eta_{n-1})}_{h_{n-1}} \otimes  \cr
   {}  &  \hskip-5pt   \otimes \widehat{L}^-_{+\alpha_1} \zeta^{\,
(1)}_{+\alpha_1} \otimes \cdots \otimes \widehat{L}^-_{+\alpha_1}
\zeta^{(\lambda^+_{+\alpha_1})}_{+\alpha_1} \otimes \cdots \otimes
\widehat{L}^-_{+\alpha_N} \zeta^{\,(1)}_{+\alpha_N} \otimes \cdots \otimes
\widehat{L}^-_{+\alpha_N} \zeta^{(\lambda^+_{+\alpha_N})}_{+\alpha_N}
\otimes  \cr
   {}  &  \hskip13pt  \otimes \psi_1 \otimes \cdots \otimes
\psi_{d - |\Cal{M}|} \;\; + \;\; T  \cr }  $$
where  $ \, \alpha_1 \preceq \alpha_2 \preceq \cdots \preceq
\alpha_N \, $  (with  $ \, N = {n \choose 2} $)  are the
positive roots of  $ \gergln \, $,  \, each one of the
$ \zeta^{(k)}_{-\alpha_r} $'s,  the  $ \theta^{(s)}_{h_i} $'s,  the
$ \zeta^{(\ell)}_{+\alpha_r} $'s  and the  $ \psi_p $'s  is a suitable
monomial in the  $ g_j^{\pm 1} $'s,  and finally  $ T $  is a sum of
homogeneous terms whose degrees are different from the degree of
the previous summand.  From this and  $ \, \epsilon \big(
\widehat{L}^\pm_{i,j} \big) = 0 = \epsilon(\ell_k) \, $  (for
all  $ k $  and all  $ \, i \not= j \, $)  we argue
  $$  \hbox{ $ \eqalign{
   \delta_d\big(\Cal{M}\big) \,  &  = \,
\widehat{L}^+_{-\alpha_1} \zeta^{\,(1)}_{-\alpha_1} \otimes \cdots \otimes
\widehat{L}^+_{-\alpha_1} \zeta^{(\lambda^+_{\alpha_1})}_{-\alpha_1}
\otimes \cdots \otimes \widehat{L}^+_{-\alpha_N} \zeta^{\,(1)}_{-\alpha_N}
\otimes \cdots \otimes \widehat{L}^+_{-\alpha_N}
\zeta^{(\lambda^+_{-\alpha_N})}_{-\alpha_N} \otimes   \hfill  \cr
   {} \otimes h_1 \,  &  \theta^{(1)}_{h_1} \otimes \cdots
\otimes h_1 \theta^{(\eta_1)}_{h_1} \otimes \cdots \otimes
h_{n-1} \theta^{(1)}_{h_{n-1}} \otimes \cdots \otimes h_{n-1}
\theta^{(\eta_{n-1})}_{h_{n-1}} \otimes  \cr
   {}  &  \hskip-8pt   \otimes \widehat{L}^-_{+\alpha_1}
\zeta^{\,(1)}_{+\alpha_1} \otimes \cdots \otimes \widehat{L}^-_{+\alpha_1}
\zeta^{(\lambda^+_{+\alpha_1})}_{+\alpha_1} \otimes \cdots \otimes
\widehat{L}^-_{+\alpha_N} \zeta^{\,(1)}_{+\alpha_N} \otimes \cdots \otimes
\widehat{L}^-_{+\alpha_N} \zeta^{(\lambda^+_{+\alpha_N})}_{+\alpha_N}
\otimes  \cr
   {}  &  \hskip7pt  \otimes (\psi_1 - 1) \otimes \cdots \otimes
(\psi_{d - |\Cal{M}|} - 1) \;\; + \;\, P  \cr } $ }   (6.1)  $$ 
where  $ \, P := {(\id - \epsilon)}^{\otimes d}\big(T\big) \, $  is
again a sum of homogeneous terms whose degrees are different from that
of the previous summand (which is homogeneous too).  In the latter
each tensor factor belongs to  $ \, U_\h(\gergln) \setminus \h \,
U_\h(\gergln) \, $,  \, whilst  $ \, (\psi_k - 1) \in \h \, U_\h
(\gergln) \setminus \h^2 \, U_\h(\gergln) \, $  for all  $ k \, $:  \,
the outcome is  $ \; \delta_d\big(\Cal{M}\big) \in \h^{d - |\Cal{M}|}
\, U_\h(\gergln) \setminus \h^{d - |\Cal{M}| + 1} \, U_\h(\gergln) \; $
for all  $ \, d \geq \big|\Cal{M}\big| \, $,  \; whence
  $$  \widetilde{\Cal{M}} \, := \, \h^{|\Cal{M}|} \, \Cal{M} \,
\in \, {U_\h(\gergln)}' \setminus \h \, {U_\h(\gergln)}'  \qquad
\forall \;\; \Cal{M} \in \Bbb{S} \; .  $$
From this we eventually get  $ \; \widetilde{\Bbb{S}} := \Big\{
\widetilde{\Cal{M}} \;\Big|\, \Cal{M} \in \Bbb{S} \,\Big\} \subseteq
{U_\h(\gergln)}' \, $,  \; thus also the  $ \kh $--\,span  of  $ \,
\widetilde{\Bbb{S}} \, $  is contained in  $ {U_\h(\gergln)}' \, $.
In fact, the previous analysis also allows to revert this last result,
thus proving the following
%
%
 \eject

   {\sl  $ \underline{\hbox{\sl Claim}} \, $:}  {\it  $ \;
\widetilde{\Bbb{S}} $  is a topological\/  $ \kh $--basis
of  $ \, {U_\h(\gergln)}' \, $.}

 \vskip4pt

   Indeed, let  $ \, \eta \in {U_\h(\gergln)}' \, $  and take an
expansion  $ \; \eta = \sum_{\Cal{M} \in \Bbb{S}} c_{\Cal{M}} \,
\Cal{M} \; $  of  $ \eta $  of minimal length as a linear combination
over  $ \kh $  of elements of  $ \Bbb{S} \, $.  Let's call  $ \,
\Cal{M}^{d,\otimes} \, $  the first summand in right-hand-side
of (6.1): then our analysis gives
 \vskip-15pt
  $$  \delta_d(\eta) \, = \, {\textstyle \sum_{\Cal{M} \in \Bbb{S}}}
\, c_{\Cal{M}} \, \delta_d\big(\Cal{M}\big) \, = \, {\textstyle
\sum_{\Cal{M} \in \Bbb{S}}} \, c_{\Cal{M}} \, \big( \Cal{M}^{d,\otimes}
+ P \big) \, = \, {\textstyle \sum_{|\Cal{M}| = \mu_+}} \, c_{\Cal{M}}
\, \Cal{M}^{d,\otimes} \, + \; R_-  $$
 \vskip-5pt
\noindent
 where  $ \, \mu_+ := \max \big\{ |\Cal{M}| \big\}_{\Cal{M} \in \Bbb{S}}
\; $  and  $ R_- $  is a sum of homogeneous terms whose degrees are
different from the degrees of any summand in  $ \, {\textstyle
\sum_{|\Cal{M}| = \mu_+}} c_{\Cal{M}} \, \Cal{M}^{d,\otimes}
\, $.  Therefore  $ \, \delta_d(\eta) \in \h^d \, U_\h(\gergln)
\, $  (as  $ \, \eta \in {U_\h(\gergln)}' \, $)  forces also  $ \,
{\textstyle \sum_{|\Cal{M}| = \mu_+}} \! c_{\Cal{M}} \, \Cal{M}^{d,
\otimes} \in \! \h^d \, U_\h(\gergln) \, $.  Again by a simple
degree argument we get  $ \; \sum_{\Sb  |\Cal{M}| = \mu_+  \\
                             \partial(\Cal{M}) = \beta  \endSb}
\hskip-5pt c_{\Cal{M}} \, \Cal{M}^{d,\otimes} \in \! \h^d
\, U_\h(\gergln) \; $  for all  $ \, \beta \in Q \, $.  Using linear
independence of monomials in the  $ \widehat{L}^\pm_\alpha $'s  with
different exponents (consequence of the quantum PBW theorem) we get
also  $ \sum\limits_{\Cal{M} \in S_{\mu_+}} \hskip-7pt c_{\Cal{M}} \,
\Cal{M}^{d,\otimes} \in \h^d \, U_\h(\gergln) \; $  where  $ S_{\mu_+} $
is the set of all monomials  $ \Cal{M} $  with  $ \, |\Cal{M}| = \mu_+
\, $  and  {\sl fixed\/}  exponents  $ \lambda^\pm_\alpha \, $.  Again
by quantum PBW, this happens if and only if  $ \sum\limits_{\Cal{M}
\in S_{\mu_+}} \hskip-7pt  c_{\Cal{M}} \, \Cal{M} \in \h^{\,\mu_+} \,
U_\h(\gergln) \, $,  \; which in turn implies  $ \, c_{\Cal{M}} \in
\h^{\,\mu_+} \, \kh \, $  for all  $ \Cal{M} $  involved; so this last
sum can be written as  $ \; \eta_+ \! = \hskip-7pt \sum\limits_{\Cal{M}
\in S_{\mu_+}} \hskip-7pt c_{\Cal{M}} \, \Cal{M} = \hskip-7pt
\sum\limits_{\Cal{M} \in S_{\mu_+}} \hskip-7pt \widehat{c}_{\Cal{M}}
\, \widetilde{\Cal{M}} \, $,  \; which belongs to the topological
$ \kh $--\,span  of  $ \widetilde{\Bbb{S}} \, $,  with  $ \,
\widehat{c}_{\Cal{M}} := \h^{-\mu_+} c_{\Cal{M}} \in \kh \, $.
But then also  $ \; \eta' := \eta - \eta_+ \in {U_\h(\gergln)}' \, $,
\; and  $ \eta' $  has less non-zero coefficients in its expansion
w.r.t.~the topological  $ \kh $--basis  $ \Bbb{S} \, $.  Iterating
this argument, we eventually find that  $ \eta $  belongs to the
topological  $ \kh $--\,span  of  $ \, \widetilde{\Bbb{S}} \, $,
\, q.e.d.
                                            \par
   Note that each  $ \, \widetilde{M} \in \widetilde{\Bbb{S}} \, $  is a
monomial in the elements  $ \, \widetilde{\ell}_k := \h \, \ell_k \, $ 
and the  $ \, \h \, \widehat{L}^\pm_{\mp \, \alpha} = \h \, \big( q - q^{-1}
\big)^{-1} L^\pm_{\mp\,\alpha} \, $,  \, hence  {\sl these are topological
algebra generators for  $ {U_\h(\gergln)}' \, $}.  Furthermore, since
$ \, \h^{-1} \big( q - q^{-1} \big) \, $  is an invertible element of
$ \kh \, $,  \, we have also that  $ {U_\h(\gergln)}' $  {\sl is generated,
as a unital  $ \kh $--\,algebra,  by the  $ L^\pm_{i,j} $'s  and the
$ \widetilde{\ell}_k $'s  (for all  $ i $,  $ j $,  $ k \, $)}.
                                            \par
   In the semiclassical limit  $ \; {U_\h(\gergln)}' \Big|_{\h=0}
\! \cong F\big[\big[GL_n^{\,*}\big]\big] = F \big[ \big[
B^{\scriptscriptstyle G}_+ \star B^{\scriptscriptstyle G}_-
\big] \big] = F \big[ \big[ \frak{b}^{\scriptscriptstyle G}_+ \star
\frak{b}^{\scriptscriptstyle G}_- \big] \big] \, $,  \; the above generators
specialize to matrix coefficients onto  $ \, \frak{b}^{\scriptscriptstyle
G}_+ \star \frak{b}^{\scriptscriptstyle G}_- \, $;  \, hereafter
$ B^{\scriptscriptstyle G}_\pm $  is the Borel subgroup in  $ GL_n $  of
upper/lower triangular matrices and  $ \, \frak{b}^{\scriptscriptstyle
G}_\pm := \text{\it Lie} \big( B^{\scriptscriptstyle G}_\pm \big) \, $, 
\, so  $ \, B^{\scriptscriptstyle G}_+ \star B^{\scriptscriptstyle G}_-
\, $  is the Poisson group dual to  $ GL_n^{\,*} $,  \, and we identify 
$ \, B^{\scriptscriptstyle G}_+ \star B^{\scriptscriptstyle G}_- \cong
\frak{b}^{\scriptscriptstyle G}_+ \star \frak{b}^{\scriptscriptstyle G}_-
\, $  (everything is very similar to the case of  $ SL_n $).  Namely, for
every  $ \, i < j \, $  the coset modulo  $ \, \h \, {U_\h(\gergln)}' \, $ 
of each $ L^+_{i,j} $  is the matrix coefficient  $ e_{i,j} $  onto  $ \,
\big(\frak{b}^{\scriptscriptstyle G}_+ \, , \, \boldkey{0} \big) \cong
\frak{b}^{\scriptscriptstyle G}_+ \, $,  \, and the coset of each
$ L^-_{j,i} $  is the matrix coefficient  $ e_{j,i} $  onto  $ \,
\big( \boldkey{0} \, , \frak{b}^{\scriptscriptstyle G}_- \big) \cong
\frak{b}^{\scriptscriptstyle G}_- \, $;  \, also, for each  $ k $  the
coset of  $ \widetilde{\ell}_k $  modulo  $ \, \h \, {U_\h(\gergln)}'
\, $  is  $ \, e_{k,k}\Big|_{B^{\scriptscriptstyle G}_+} \! = \,
e_{k,k}^{\;-1}\Big|_{B^{\scriptscriptstyle G}_-} \; $.  Finally,
as  $ \, L^\pm_{k,k} \! =: \! g_k^{\pm 1} \! := \exp \! \big( \h \,
\ell_k \big) \! = \exp \! \big(\,\widetilde{\ell}_k\big) \, $  the
same kind of relation occurs between the cosets modulo  $ \, \h \,
{U_\h(\gergln)}' \, $  of  $ L^\pm_{k,k} $  and of  $ \widetilde{\ell}_k
\, $,  for all  $ k \, $.
                                            \par
   As for  $ {U_\h(\gersln)}' $,  for all  $ \, i < j \, $  we have that
$ \; \widetilde{F}_{j,i} := \big( q - q^{-1} \big) \, F_{j,i} = g_i^{-1}
\, L^+_{i,j} \; $  and  $ \; \widetilde{E}_{i,j} := - \, \big( q - q^{-1}
\big) \, E_{i,j} = - L^-_{j,i} \, g_i^{+1} \; $  belong to  $ \,
{U_\h(\gergln)}' \, \bigcap \, U_\h(\gersln) = {U_\h(\gersln)}' \, $,
\, as well as  $ \, \widetilde{h}_k := \h \, (\ell_k - \ell_{k+1}) =
\widetilde{\ell}_k - \widetilde{\ell}_{k+1} \, $  (for all  $ k \, $).
Indeed, with the same analysis as above   --- up to the obvious, minimal
changes ---   one proves also that  $ {U_\h(\gersln)}' $  {\sl is
generated, as a topological unital  $ \kh $--\,algebra,
     \hbox{by the  $ \widetilde{F}_{j,i} $'s,  the
$ \widetilde{E}_{i,j} $'s  (for all  $ \, i < j \, $)
and the  $ \widetilde{h}_k $'s  (for all  $ k \, $)}.}
%
%
%
 \eject
   In addition,  $ {U_\h(\gersln)}' $  has  $ \kh $--basis  the set of
rescaled PBW-like monomials (in the above generators) analogue to the
set  $ \, \widetilde{\Bbb{S}} \, $  considered above which is a basis
for  $ {U_\h(\gergln)}' \, $.
                                        \par
   Finally, under specialization  $ \; {U_\h(\gersln)}'\Big|_{\h=0}
\! \cong F\big[\big[SL_n^{\,*}\big]\big] = F\big[\big[ B_+ \star B_-
\big]\big] = F\big[\big[\frak{b}_+ \star \frak{b}_-\big]\big] \, $ 
the above generators specialize as  $ \, \widetilde{F}_{j,i}\big|_{\h=0}
\! = e_{i,i}^{\;-1} \, e_{i,j} \big|_{\frak{b}_+} \, $,  $ \,
\widetilde{E}_{i,j}\big|_{\h=0} \! = e_{j,i} \; e_{i,i}^{\;+1}
\big|_{\frak{b}_-} \, $  (for all  $ \, i < j \, $)  and  $ \,
\widetilde{h}_k\big|_{\h=0} \! = e_{k,k} \big|_{\frak{b}_+} -
e_{k+1,k+1} \big|_{\frak{b}_+} \, $  (for all  $ \, k = 1,
\dots, n-1 \, $).

 \vskip7pt

   {\bf 6.9 Quantum Stokes matrices:  $ \, \boldkey{n} \boldkey{=}
\boldkey{3} \, $.} \, According to the general recipe in \S 6.7, the
generators of  $ \, \Cal{H} = U_\h(\gerso_3) \, $  are
  $$  \displaylines{
   K_{1,2} \, = \, q^2 \, \big( F_1 - q \, T_1^{-1} E_1 \big) \;\; ,
\qquad \qquad  K_{2,3} \, = \, q \, \big( F_2 - q \, T_2^{-1} E_2 \big)
\cr
   K_{1,3} \, = \, q^2 \, \big( F_{3,1} - \big( q - q^{-1} \big) \,
F_2 \, T_1^{-1} E_1 - T_1^{-1} T_2^{-1} E_{1,3} \big)  \cr }  $$
(cf.~\S 6.7) where  $ \, T_s^{\pm 1} := t_s^{\pm 1} \, $
($ s = 1, 2 $).  From this one can directly prove that
  $$  \big[ K_{1,2} \, , K_{2,3} \big]_q \, =
\, - \, q^2 \, K_{1,3} \;\; .   \eqno (6.2)  $$
Using the relations between the elements  $ \theta_j $  in [No,
\S 2.4] --- namely, formulas (2.23) therein ---   and remarking
that  $ \; K_{1,2} = q \, \theta_1 \, $,  $ \, K_{2,3} = \theta_2
\, $,  \; one can derive also
  $$  \big[ K_{1,3} \, , K_{1,2} \big]_q \, = \, -
\, q^3 \, K_{2,3} \;\; ,  \qquad
\big[ K_{2,3} \, , K_{1,3} \big]_q \, = \, - \,
q \, K_{1,2} \;\; .   \eqno (6.3)  $$
Indeed, the case  $ \, n = 3 \, $  is especially interesting because,
using renormalized generators  $ \widetilde{K}_{1,2} := q^{-5/2}
K_{1,2} \, $,  $ \widetilde{K}_{1,3} := q^{-4/2} K_{1,3} \, $  and
$ \widetilde{K}_{2,3} := q^{-3/2} K_{2,3} \, $  one has for  $ \,
U_\h(\gerso_3) \, $  a cyclically invariant presentation (see [HKP]
and references therein, and  Remark 6.11{\it (b)\/}  too).  However,
this special feature has no general counterpart for  $ \, n \not= 3 \, $.
                                            \par
   The following PBW-like theorem holds for  $ U_\h(\gerso_3) \, $,
\, as a direct consequence of definitions and formulas (6.2--6.3):

\vskip4pt

   {\sl  $ \underline{\hbox{\sl Claim}} \, $:}  {\it  $ \;
U_\h(\gerso_3) $  is a topologically free\/  $ \kh $--\,module,
with topological  $ \kh $--\,basis  the set of ordered monomials\/
$ \; \big\{ K_{1,2}^{\,a} \, K_{1,3}^{\,b} \, K_{2,3}^{\,c} \,\big|\,
a, b, c \in \N \big\} \; $.  A similar basis is the one with
$ \widetilde{K}_{i,j} $  instead of  $ K_{i,j} $  everywhere.}

\vskip7pt

\proclaim{Theorem 6.10}  $ \, F_\h\big[\big[U^+_3\big]\big] :=
{U_\h(\gerso_3)}^{\!\Lsh} \, $  is the topological,  $ \h $--adically
complete, unital\/  $ \kh $--\,algebra  with generators
  $$  k_{1,2} \, := \, q^{-2} \, \big( q - q^{-1} \big) \, K_{1,2} \;\; ,
\quad  k_{2,3} \, := \, q^{-1} \, \big( q - q^{-1} \big) \, K_{2,3} \;\; ,
\quad  k_{1,3} \, := \, q^{-2} \, \big( q - q^{-1} \big) \, K_{1,3}  $$
and relations
  $$  \hbox{ $ \eqalign{
   k_{1,2} \, k_{2,3} \;  &  = \; q \, k_{2,3} \, k_{1,2} \, -
\, q \, \big( q - q^{-1} \big) \, k_{1,3}  \cr
   k_{2,3} \, k_{1,3} \;  &  = \; q \, k_{1,3} \, k_{2,3} \, -
\, \big( q - q^{-1} \big) \, k_{1,2}  \cr
   k_{1,3} \, k_{1,2} \;  &  = \; q \, k_{1,2} \, k_{1,3} \, -
\, \big( q - q^{-1} \big) \, k_{2,3}  \cr } $ }   \eqno (6.4)  $$
with the right coideal structure given by
  $$  \displaylines{
   \Delta\big(k_{1,2}\big) \, = \, 1 \otimes k_{1,2} \, + \,
k_{1,2} \otimes t_1^{-1} \;\; ,  \qquad \quad
 \Delta\big(k_{2,3}\big) \, = \, 1 \otimes k_{2,3} \, + \,
k_{2,3} \otimes t_2^{-1}  \cr
   \Delta\big(k_{1,3}\big) \, = \, 1 \otimes k_{1,3} \, + \,
k_{1,3} \otimes t_1^{-1} t_2^{-1} + \, \big( q \!-\! q^{-1} \big) \,
k_{1,2} \otimes f_2 \, t_1^{-1} - \, q^{-1} \big( q \!-\! q^{-1} \big)
\, k_{2,3} \otimes t_1^{-1} t_2^{-1} e_1 \, .  \cr }  $$
Moreover,  $ \, F_\h\big[\big[U^+_3\big]\big] := {U_\h(\gerso_3)}^{\!\Lsh}
\, $  is a free\/  $ \kh $--\,module,  a  $ \kh $--\,basis  being the
set of ordered monomials\/  $ \, \Bbb{B}_3 := \big\{ k_{1,2}^{\,a} \,
k_{1,3}^{\,b} \, k_{2,3}^{\,c} \,\big|\, a, b, c \in \N \,\big\} \, $.
\endproclaim
 \eject

\demo{Proof} The relations (6.4) among the  $ k_{i,j} $'s  clearly
spring out of formulas (6.2)--(6.3), whilst the formulas for the
right coideal structure directly come out of the very definitions. 
The key point of the proof instead is to show that these elements
do generate  $ {U_\h(\gerso_3)}^{\!\Lsh} \, $.
                                         \par
   From the above formulas for  $ \Delta \, $,  a straightforward
computation proves that  ($ \, \forall \; d \in \N \, $)
  $$  \displaylines{
   \delta_d\big(k_{1,2}\big) \, = \, k_{1,2} \otimes
\big( t_1^{-1} - 1 \big)^{\otimes (d-1)} \;\; ,  \qquad \quad
   \delta_d\big(k_{2,3}\big) \, = \, k_{2,3} \otimes
\big( t_2^{-1} - 1 \big)^{\otimes (d-1)}  \cr
   \quad  \delta_d\big(k_{1,3}\big) \, = \, k_{1,3} \otimes
\big( t_1^{-1} t_2^{-1} - 1 \big)^{\otimes (d-1)} \, +   \hfill  \cr
   {} \hfill   + \, {\textstyle \sum\nolimits_{r + s = d - 2}} \,
\big( q - q^{-1} \big) \, k_{1,2} \otimes \big( t_1^{-1} - 1
\big)^{\otimes r} \otimes f_2 \, t_1^{-1} \otimes \big( t_1^{-1}
t_2^{-1} - 1 \big)^{\otimes s} \, +  \qquad \quad  \cr
   {} \hfill   + \, {\textstyle \sum\nolimits_{r + s = d - 2}}
\; q^{-1} \big( q - q^{-1} \big) \, k_{2,3} \otimes \big( t_2^{-1}
- 1 \big)^{\otimes r} \otimes t_1^{-1} t_2^{-1} e_1 \otimes
\big( t_1^{-1} t_2^{-1} - 1 \big)^{\otimes s}  \cr }  $$
As  $ \; k_{i,j} \, , \, \big( t_h^{-1} - 1 \big) \, , \, \big( t_1^{-1}
t_2^{-1} - 1 \big) \in \h \, U_\h(\gerso_3) \setminus \h^2 \, U_\h
(\gerso_3) \, $,  \, we have  $ \; k_{1,2} \, , k_{2,3} \, , k_{1,3}
\in {U_\h(\gerso_3)}^{\!\Lsh} \setminus \h \, {U_\h(\gerso_3)}^{\!\Lsh}
\, $,  \; so the subalgebra generated by these elements lies in
$ {U_\h(\gerso_3)}^{\!\Lsh} \, $.
                                       \par
   We shall now prove that  $ \Bbb{B}_3 $  is a topological
$ \kh $--\,basis  of  $ {U_\h(\gerso_3)}^{\!\Lsh} \, $;  \,
this in turn will imply that this algebra is generated by  $ \,
k_{1,2} \, $,  $ \, k_{2,3} \, $  and  $ \, k_{1,3} \, $.  First,
the  $ \underline{\hbox{\sl Claim\/}} $  in \S 6.9 implies that
$ \Bbb{B}_3 $  is a linearly independent set inside  $ {U_\h
(\gerso_3)}^{\!\Lsh} \, $;  \, then now we prove that it spans
$ {U_\h(\gerso_3)}^{\!\Lsh} $  over  $ \kh \, $.  The formulas
for  $ \Delta $  on the  $ k_{i,j} $'s  give also, for all
$ \, d \in \N \, $,
  $$  \displaylines{
   \Delta^d\big(K_{1,2}\big) \, = \, {\textstyle
\sum\limits_{r + s = \, d - 1}} \hskip-3pt 1^{\otimes r} \otimes
K_{1,2} \otimes \big( t_1^{-1} \big)^{\otimes s} \, ,  \quad
   \Delta^d\big(K_{2,3}\big) \, = \, {\textstyle
\sum\limits_{r + s = \, d - 1}} \hskip-3pt 1^{\otimes r} \otimes
K_{2,3} \otimes \big( t_2^{-1} \big)^{\otimes s}  \cr
   \quad   \Delta^d\big(K_{1,3}\big) \, = \, {\textstyle
\sum\limits_{r + s = \, d - 1}} \hskip-3pt 1^{\otimes r} \otimes
K_{1,3} \otimes \big( t_1^{-1} t_2^{-1} \big)^{\otimes s}
\, +   \hfill  \cr
   + \, {\textstyle \sum\limits_{r + p + s = \, d - 2}}
\hskip-3pt 1^{\otimes r} \otimes K_{1,2} \otimes
\big( t_1^{-1} \big)^{\otimes p} \otimes A \otimes
\big( t_1^{-1} t_2^{-1} \big)^{\otimes s} \, +  \cr
   \hfill   + \, {\textstyle
\sum\limits_{r + p + s = \, d - 2}} \hskip-3pt 1^{\otimes r} \otimes
K_{2,3} \otimes \big( t_2^{-1} \big)^{\otimes p} \otimes B
\otimes \big( t_1^{-1} t_2^{-1} \big)^{\otimes s}  \cr
   \hbox{with}  \,  A := L^+_{2,3} \, g_1^{-1} \! = \big( q - q^{-1} \big)
f_2 \, t_1^{-1} ,  \,  B := q^{-1} g_3 \, L^-_{2,1} \! = - \, q^{-1} \!
\big( q - q^{-1} \big) \, t_1^{-1} t_2^{-1} e_1 \! \in {U_\h(\gersln)}'.
\cr }  $$
In particular, this implies that  $ \; \delta_{a+2\,b+c}\big(K_{1,2}^{\,a}
\, K_{1,3}^{\,b} \, K_{2,3}^{\,c} \big) = \sum_{i \in I} C_{i,1} \otimes
C_{i,2} \otimes \cdots \otimes C_{i,a+2\,b+c} \; $  (for some index set
$ I \, $)  where each tensor factor  $ C_{i,j} $  is a product of type
  $$  C_{i,j} \, = \, t_1^{-n_1} t_2^{-\nu_1} \cdot D_1 \cdot t_1^{-n_2}
t_2^{-\nu_2} \cdot D_2 \cdot \cdots \cdot t_1^{-n_{k-1}} t_2^{-\nu_{k-1}}
\cdot D_{k-1} \cdot t_1^{-n_k} t_2^{-\nu_k}   \eqno (\, k \in \N_+)  $$
with  $ \, n_s, \nu_s \in \N \, $  and  $ \, D_s \in \big\{ K_{1,2} \, ,
K_{1,3} \, , K_{2,3} \, , A \, , B \, \big\} \bigcup \big\{\! \big(
t_1^{-\tau_1} t_2^{-\tau_2} - 1 \big) \,\big|\, \tau_1, \tau_2 \in
\N_+ \big\} \, $.  In particular   --- cf.~also (6.4) ---   there
is a first summand of type
  $$  \Phi_1^{a,b,c} \!\! = \! \bigg( {\textstyle
\bigotimes\limits_{p=0}^{a-1}} t_1^{-p} K_{1,2} \!\bigg) \!
\otimes \! \left( {\textstyle \bigotimes\limits_{r=0}^{b-1}}
t_1^{-(a+r)} t_2^{-r} \!\right) \! \otimes
 \left( {\textstyle \bigotimes\limits_{s=0}^{c-1}}
t_1^{-(a+b)} t_2^{-(b+s)} K_{2,3} \!\right) \! \otimes
 \big( t_1^{-(a+b)} t_2^{-(b+c-1)} \! - 1 \big)^{\! \otimes b}  $$
   \indent   Define the  {\sl length\/}  of  $ \; K_{1,2}^{\,a} \,
K_{1,3}^{\,b} \, K_{2,3}^{\,c} \in \Bbb{B}_3 \; $  as  $ \; l \big(\,
K_{1,2}^{\,a} \, K_{1,3}^{\,b} \, K_{2,3}^{\,c} \big) \, := \, a + 2 \,
b + c \, $,  \; and let  $ \, \Cal{H}_n \, $  be the  $ \kh $--\,span
of all monomials in  $ \Bbb{B}_3 $  of length at most  $ n \, $.  This
defines an algebra filtration  $ \, \{\Cal{H}_n\}_{n \in \N} \, $  of
$ U_\h(\gerso_3) \, $;  \, the formulas for the coproduct of the
$ k_{i,j} $'s  show that this is a  {\sl comodule algebra filtration},
i.e.~an algebra filtration such that  $ \, \Delta \big( \Cal{H}_n \big)
\subseteq \Cal{H}_n \otimes U_\h(\gersl_3) \, $  for all  $ n \, $.
A similar filtration is also induced onto each tensor power
$ {U_\h(\gerso_3)}^{\otimes\, l} $  ($ \, l \in \N \, $).
                                          \par
  Any  $ \, \eta \in {U_\h(\gerso_3)}^{\!\Lsh} \, $  expands uniquely
as  $ \; \eta = \sum_{a,b,c \in \N} \chi_{a,b,c} \, K_{1,2}^{\,a} \,
K_{1,3}^{\,b} \, K_{2,3}^{\,c} \; $  for some  $ \, \chi_{a,b,c} \in
\kh \, $,  \, by the  $ \underline{\hbox{\sl Claim\/}} $  of \S 6.9.
Set  $ \, \mu := \min \big\{ a + 2 \, b + c \,\big|\, \chi_{a,b,c}
\not= 0 \big\} \, $,  \, and look at  $ \; \delta_\mu(\eta) =
\sum_{a,b,c \in \N} \chi_{a,b,c} \cdot \delta_\mu \big( K_{1,2}^{\,a}
\, K_{1,3}^{\,b} \, K_{2,3}^{\,c} \big) \in \h^{\,\mu} \, U_\h(\gerso_3)
\otimes {U_\h(\gersl_3)}^{\otimes (\mu - 1)} \, $.
             By degree\break
 \eject
\noindent
 arguments   --- w.r.t.~the filtration  $ \{\Cal{H}_n\}_{n \in \N} $
of  $ U_\h(\gerso_3) $  given above ---   we see that  $ \,
\delta_\mu(\eta) \in \h^{\,\mu} \, U_\h(\gerso_3) \otimes
{U_\h(\gersl_3)}^{\otimes (\mu - 1)} \, $  forces also
  $$  {\textstyle \sum\limits_{a+2\,b+c=\mu}} \chi_{a,b,c} \cdot
\delta_\mu \big( K_{1,2}^{\,a} \, K_{1,3}^{\,b} \, K_{2,3}^{\,c} \big)
\, \in \, \h^{\,\mu} \, U_\h(\gerso_3) \otimes {U_\h(\gersl_3)}^{\otimes
(\mu - 1)} \; .   \eqno (6.5)  $$
   \indent   By the analysis above, each  $ \, \delta_\mu \big(
K_{1,2}^{\,a} \, K_{1,3}^{\,b} \, K_{2,3}^{\,c} \big) \, $  in (6.5)
is equal to  $ \, \Phi_1^{a,b,c} \, $  (defined above) plus other
terms which are linearly independent of  $ \Phi^{a,b,c} $  modulo  $ \h
\, {U_\h(\gersln)}^{\otimes \mu} \, $.  Furthermore, all these
$ \Phi_1^{a,b,c} $'s,  for different triples  $ \, (a,b,c) \in \N^3
\, $,  \, are linearly independent inside  $ {U_\h(\gersln)}^{\otimes
\mu} \, $,  \, by construction.  As an outcome, we have that (6.5)
implies
  $$  \chi_{a,b,c} \cdot \Phi_1^{a,b,c} \, \in \, \h^{\,\mu} \,
U_\h(\gerso_3) \otimes {U_\h(\gersl_3)}^{\otimes (\mu - 1)}
\eqno \forall \;\; a + 2 \, b + c = \mu \; .  $$
Since  $ \, \Phi_1^{a,b,c} \in \h^b \, U_\h(\gerso_3) \otimes
{U_\h(\gersl_3)}^{\otimes (\mu - 1)} \, $  by construction, we
argue  $ \, \chi_{a,b,c} \in \h^{a+b+c} \, \kh \, $  for all
$ \, a + 2 \, b + c = \mu \, $,  \, so that
  $$  \chi_{a,b,c} \, K_{1,2}^{\,a} \, K_{1,3}^{\,b} \, K_{2,3}^{\,c}
\, \in \, \kh \cdot k_{1,2}^{\,a} \, k_{1,3}^{\,b} \, k_{2,3}^{\,c} \,
\subseteq \, \kh\hbox{--\,span of} \ \Bbb{B}_3   \eqno \forall \;\;
a + 2 \, b + c = \mu \; .  $$
But then  $ \, \eta_- := \sum\limits_{a+2\,b+c=\mu} \hskip-7pt
\chi_{a,b,c} \cdot K_{1,2}^{\,a} \, K_{1,3}^{\,b} \, K_{2,3}^{\,c}
\, \in \, {U_\h(\gerso_3)}^{\!\Lsh} \, $  by our previous results,
hence also
  $$  \eta_> \, := \, \eta - \eta_- \, = \,
{\textstyle \sum\limits_{a + 2 \, b + c > \mu}} \hskip-5pt
\chi_{a,b,c} \cdot K_{1,2}^{\,a} \, K_{1,3}^{\,b} \, K_{2,3}^{\,c}
\, \in \, {U_\h(\gerso_3)}^{\!\Lsh} \; .  $$
Now we can apply the same arguments to  $ \eta_< $  instead of
$ \eta \, $:  \, iterating this procedure (involving monomials in
the  $ K_{i,j} $'s  whose length grows up), we eventually find
that  $ \eta $  belongs to the topological  $ \kh $--\,span
of  $ \Bbb{B}_3 \, $,  \, q.e.d.   \qed
\enddemo

\vskip7pt

   {\bf 6.11 Remarks:} {\it (a)} \, in \S 6.8 we saw that
$ {U_\h(\gersln)}' $  is generated by the  $ L $--operators,
hence its semiclassical limit  $ F[[G^*]] $  is generated by their
cosets, which are simply half the matrix coefficients generating
$ F[[G^*]] $  (see \S 6.1).  Then by the very construction and our
concrete description of  $ {U_\h(\gerso_3)}^{\!\Lsh} \, $  we get that
the generators  $ k_{i,j} $  specialize, in  $ \, {U_\h(\gerso_3)}^{\!
\Lsh}{\Big|}_{\h=0} = F\big[\big[U^+_3\big]\big] \, $,  \, right
to the generators of  $ \, F\big[\big[U^+_3\big]\big] \, $  (cf.~\S
6.1).  In particular, the corresponding limit Poisson bracket can
therefore be verified to be equal to that in [Ug] and in [Xu] (the
latter taken from [Du]), up to normalizations: e.g., the isomorphism
between our presentation of  $ F\big[\big[U^+_3\big]\big] $  and
Xu's one is given by
  $$  k_{1,2}{\big|}_{\h=0} \mapsto z \;\; ,  \qquad
k_{1,3}{\big|}_{\h=0} \mapsto y \;\; ,  \qquad
k_{2,3}{\big|}_{\h=0} \mapsto x  $$
(notation of [Xu], \S 1, formula (2)), and this is easily
       \hbox{seen to preserve the Poisson bracket.}
                                                 \par
   {\it (b)} \, the claim and proof of Theorem 6.10 show that one
could take as generators for  $ {U_\h(\gerso_3)}^{\!\Lsh} \, $  simply
the  $ \, \big( q - q^{-1} \big) \, K_{i,j} \, $'s.  However, our choice
of normalization (dividing out such generators by suitable powers of
$ q \, $)  lead us to better looking relations, such as (6.4).  Indeed,
this can still be improved, taking new generators  $ \; \widetilde{k}_{1,2}
:= q^{-1/2} \, k_{1,2} = \big( q - q^{-1} \big) \widetilde{K}_{1,2} \; $,
$ \; \widetilde{k}_{1,3} := k_{1,3} = \big( q - q^{-1} \big)
\widetilde{K}_{1,3} \; $  and  $ \; \widetilde{k}_{2,3} := q^{-1/2}
\, k_{2,3} = \big( q - q^{-1} \big) \widetilde{K}_{2,3} \; $  (see
\S 6.9):these enjoy the relations
 $ \; \widetilde{k}_{1,2} \,
\widetilde{k}_{2,3} \, = \, q \, \widetilde{k}_{2,3} \,
\widetilde{k}_{1,2} \, - \, \big( q - q^{-1} \big) \,
\widetilde{k}_{1,3} \, $,  $ \; \widetilde{k}_{2,3} \,
\widetilde{k}_{1,3} \, = \, q \, \widetilde{k}_{1,3} \,
\widetilde{k}_{2,3} \, - \, \big( q - q^{-1} \big) \,
\widetilde{k}_{1,2} \, $,  $ \; \widetilde{k}_{1,3} \,
\widetilde{k}_{1,2} \, = \, q \, \widetilde{k}_{1,2} \,
\widetilde{k}_{1,3} \, - \, \big( q - q^{-1} \big) \,
\widetilde{k}_{2,3} \, $,  \;
which are  {\sl totally symmetric with respect to cyclic permutations
of the indices}. Nevertheless, this special feature   --- like for
$ U_\h(\gerso_3) $  ---   has no general counterpart for  $ \, n
\not= 3 \, $.   $ \quad \diamondsuit $
%
%
 \eject

   {\bf 6.12 The general case.} \, Let us now move to the general case
$ \, n > 3 \, $.  The generators  $ K_{i,j} $  ($ i < j $)  are defined
in \S 6.7; like in the  $ \underline{\hbox{\sl Claim\/}} $  in \S 6.9,
we have a PBW-like theorem for  $ U_\h(\gerso_n) \, $:  \, namely, the
set of all ordered monomials (w.r.t.~any fixed total order of the set
of pairs  $ \big\{ (i,j) \,\big|\, i<j \,\big\} \, $)  in the
$ K_{i,j} $'s  is a topological $ \kh $--\,basis  of
$ U_\h(\gerso_n) \, $.
                                         \par
   Straightforward computations yield
  $$  \delta_d\big(K_{i,j}\big) \, = \, {\textstyle \sum\limits_I}
\;\, K_{t_1, s_1} \otimes (\id - \epsilon) \big( L^-_{t_1, t_2}
L^+_{s_1, s_2} \big) \otimes \cdots \otimes(\id - \epsilon)
\big( L^-_{t_{d-2}, i} \, L^+_{s_{d-2}, j} \big)  $$
where the set of indices is  $ \, I = \big\{\, i \leq t_{d-2} \leq \cdots
\leq t_1 < s_1 \leq \cdots \leq s_{d-2} \leq j \,\big\} \, $;  \, it is
worth pointing out that, while the  $ L $--operators  $ L^+_{i,j} $
and  $ L^-_{i,j} $  do not belong to  $ U_\h(\gersln)$  but only
to  $ U_\h(\gergln) \, $,  the products  $ \, L^-_{t_r, t_{r+1}}
L^+_{s_r, s_{r+1}} \, $  do belong to  $ U_\h(\gersln) \, $.
From this one gets easily
  $$  \delta_d\big(K_{i,j}\big) \, \in \, \h^{d-1} U_\h(\gerso_n)
\otimes {U_\h(\gersln)}^{\otimes (d-1)}   \eqno (\, i < j \, ,
\, d \in \N)  $$
whence  $ \, k_{i,j} := \big( q - q^{-1} \big) \, K_{i,j} \in
{U_\h(\gerso_n)}^{\!\Lsh} \setminus \h \, {U_\h(\gerso_n)}^{\!\Lsh}
\, $  follows at once.
                                 \par
   Indeed, with much the same analysis as in \S\S 6.9--10 one can prove
that in fact the  $ k_{i,j} $'s  (for  $ \, i < j \, $)  form a complete
set of generators for the algebra  $ {U_\h(\gerso_n)}^{\!\Lsh} \, $,
and that the set of ordered monomials in these generators is a topological
$ \kh $--\,basis  for  $ {U_\h(\gerso_n)}^{\!\Lsh} \, $.  Finding the
relations between the  $ k_{i,j} $'s  then will provide an explicit
presentation of the algebra  $ {U_\h(\gerso_n)}^{\!\Lsh} \, $,  \,
hence a quantization  $ \, F_\h\big[\big[U^+_n\big]\big] :=
{U_\h(\gerso_n)}^{\!\Lsh} \, $  of  $ F\big[\big[U^+_n\big]\big] $
with the Poisson structure given in [Ug], the analogue of
Remark 6.11{\it (a)\/}
      \hbox{holding true in the general case too.}

\vskip1,3truecm

\centerline {\bf \S\; 7  \ Generalizations }

\vskip13pt

  {\bf 7.1 Quantum duality with half quantizations.} \, In the present
work we take from scratch the datum of a pair of mutually dual quantum
groups, namely  $ \big( \fhg \, , \, \uhg \big) \, $  (cf.~\S 2.7).  In
the proofs, this assumption is exploited to apply orthogonality arguments,
for which all these are necessary (a single quantum groups would not be
enough).
                                                \par
   However, this is only a matter of choice.  Indeed, our quantum duality
principle deals with quantum subgroups which are contained either in
$ \fhg $  or in  $ \uhg \, $,  and  {\sl we might prove every step in our
discussion using only the single quantum group which is concerned},  and
only one quantum subgroup (such as  $ \calI_\h $,  or  $ \calC _\h $,
etc.) at the time, by a direct method which use no orthogonality
arguments.  To give a sample, we re-prove part of Lemma 4.2:

\vskip2pt

  {\it  $ \underline{\hbox{\it Claim:}} $ \,  let  $ \, {\calI_\h}^{\!
\curlyvee} $  and  $ \, {\calC_\h}^{\!\!\triangledown} $  be as in
Lemma 4.2.  Then  $ \; {\calI_\h}^{\!\curlyvee} \,\coideal\; {\fhg}^\vee
\, $  and  $ \,\; {\calC_\h}^{\!\!\triangledown} \,\coleq\, {\fhg}^\vee
\, $.}

\vskip-13pt

\demo{Proof}  By definition  $ \, {\calI_\h}^{\!\curlyvee} \, $  is the
left ideal of  $ {\fhg}^\vee $ generated by  $ \, \h^{-1} \calI_\h \, $,
hence it is enough to show that  $ \, \Delta \big( {\fhg}^\vee \cdot
\h^{-1} \calI_\h \big) \, \subseteq \, {\fhg}^\vee \!\otimes {\calI_\h}^{\!
\curlyvee} + {\calI_\h}^{\!\curlyvee} \!\otimes {\fhg}^\vee \, $.  Since
$ \calI_\h $  is a coideal of  $ \fhg $  (see \S 2.6), we have  $ \,
\Delta \big( {\fhg}^\vee \cdot \h^{-1} \calI_\h \big) \, \subseteq \,
\big( {\fhg}^\vee \!\otimes {\fhg}^\vee \big) \cdot \big( \fhg \otimes
\h^{-1} \calI_\h + \h^{-1} \calI_\h \otimes \fhg \big) \, \subseteq \,
{\fhg}^\vee \!\otimes {\calI_\h}^{\!\curlyvee} + {\calI_\h}^{\!\curlyvee}
\!\otimes {\fhg}^\vee \, $,  \, q.e.d.
                                             \par
%
%
   The case of  $ \, {\calC_\h}^{\!\!\triangledown} \, $  is entirely
similar.   \qed
\enddemo

\vskip7pt

  {\bf 7.2 Quantum duality with global quantizations.} \,
In this paper we use quantum groups in the sense of Definition 2.2;
in literature, these are sometimes called  {\sl local\/}  quantizations.
Instead, one can consider  {\sl global quantizations\/}:  quantum groups
like Jimbo's, Lusztig's, etc.  The latter ones differ from the former in
two respects:
 \vskip0pt
   {\it ---1)}  \, they are standard (rather than topological)
Hopf algebras;
                                           \par
   {\it ---2)}  \, they may be defined over any ring  $ R \, $,
the r\^{o}le of  $ \h $  being played by a suitable element of
that ring (the most common example is  $ \, R = \Bbbk \big[ q,
q^{-1} \big] \, $  and  $ \, \h = q - 1 \, $).
 \vskip2pt
   The first point implies that the semiclassical limit of a quantum
group of this type is either  $ U(\gerg) \, $,  for some Lie bialgebra
$ \gerg \, $,  or  $ F[G] \, $,  the algebra of regular functions on some
Poisson algebraic group  $ G \, $.  The latter is a geometrical object of
{\sl global\/}  type, thus a quantum group specializing to it carries
richer information than a QFSHA.  The second point implies that one can
consider different specializations, namely one for each point of the
spectrum of the ground ring  $ R \, $:  so this setting is richer
from an arithmetical viewpoint.
                                               \par
   Now, the present work might be written equally well in terms of  {\sl
global quantum groups\/}  and their specializations.  The only care is to
start with algebraic Poisson groups and algebraic Poisson homogeneous
spaces, instead of formal ones.  Then one defines Drinfeld-like functors
in a perfectly similar manner; the key fact is that the quantum duality
principle has a  {\sl global version\/}  (see [Ga2]) in which the recipe
given in \S 3 to define Drinfeld-like functors do make sense, up to a few
technical details, in the global framework as well.  In addition, one can
also extend our quantum duality principle for coisotropic subgroups (and
Poisson quotients) to all closed subgroups (and all homogeneous spaces):
the outcome then is that applying the so-extended Drinfeld's functors
to any closed subgroup (or homogeneous space) one always gets a  {\sl
coisotropic\/}  subgroup (or a Poisson quotient) of the dual Poisson
group, and this is again characterized in terms of involutivity
(see [CG]).

\vskip7pt

  {\bf 7.3 $ * $--structures  and quantum duality for real subgroups and
homogeneous spaces.} \, If one is interested in quantizations of real
subgroups and real homogeneous spaces, then  $ * $--structures  must be
considered on the quantum group Hopf algebras one starts from.  It is
then possible to perform all our construction in this setting, and to
formulate and  prove a version of the QDP for  {\sl real\/}  quantum
subgroups and quantum homogeneous spaces too, both in the formal and
in the global setting; see
[CG] for details.

 \vskip1,9truecm

\vskip55pt

\enddocument
\bye